\documentclass[a4paper,12pt]{article}
\usepackage{amssymb}
\usepackage{amsmath}
\usepackage{amsthm}
\def\dj{d\kern-.30em\raise1.25ex\vbox{\hrule width .3em height .03em}}
\def\Dj{\rlap{\kern-.70em\raise0.75ex\vbox{\hrule width .3em height .03em}}} 

\newtheorem{theorem}{Theorem}[section]
\newtheorem{proposition}[theorem]{Proposition}
\newtheorem{lemma}[theorem]{Lemma}
\newtheorem{cor}[theorem]{Corollary}

\newtheorem{remark}[theorem]{Remark}

\newcommand{\A }{\mathbb{A}}

\newcommand{\ABM }{{}^\cA_\B \hspace{-.1cm}\mathcal{M}}
\newcommand{\AqM }{{}^{\Aqr}\hspace{-.1cm}\mathcal{M}}
\newcommand{\Aql }{\overrightarrow{\cA }}
\newcommand{\Aqr }{\overleftarrow{\cA }}
\newcommand{\B }{\mathcal{B}}

\newcommand{\cC}{\mathcal{C}}
\newcommand{\cD}{\mathcal{D}}
\newcommand{\cA }{\mathcal{A}}

\newcommand{\cH}{\mathcal{H}}

\newcommand{\cN}{\mathcal{N}}

\newcommand{\co }{\mathrm{co}}

\newcommand{\cQ}{\mathcal{Q}}

\newcommand{\cS}{\mathcal{S}}

\newcommand{\cqg}{\mathbb{C}_q[G]}
\newcommand{\cqgl}{\mathbb{C}_q[G/L_S]}

\newcommand{\C }{\mathbb{C}}
\newcommand{\del }{\partial}
\newcommand{\delb }{{\overline{\partial}}}
\newcommand{\dif }{\mathrm{d}}

\newcommand{\Fil }{\mathcal{F}}

\newcommand{\Gduw}[1][]{\Gamma _{\del ,\mathrm{u}}^{\wedge #1}}
\newcommand{\Gdbuw}[1][]{\Gamma _{\delb ,\mathrm{u}}^{\wedge #1}}
\newcommand{\Gdifuw}[1][]{\Gamma _{\dif ,\mathrm{u}}^{\wedge #1}}

\newcommand{\gfrak}{\mathfrak{g}}
\newcommand{\gproots}{\overline{R^+_S}}
\newcommand{\gr}{\mathrm{Gr}}
\newcommand{\hfrak}{\mathfrak{h}}

\newcommand{\hght}{\mathrm{ht}}
\newcommand{\id}{\mathrm{Id}}

\newcommand{\im}{\mathrm{Im}}
\newcommand{\Ind}{I}

\newcommand{\kfrak}{\mathfrak{k}}
\newcommand{\kopr }{\varDelta }
\newcommand{\kow }{\varDelta }
\newcommand{\K }{K}
\newcommand{\Kkow }{\varDelta _\K }
\newcommand{\lact }{\triangleright}

\newcommand{\lfrak}{\mathfrak{l}}

\newcommand{\Lin }{\mathrm{Lin}}

\newcommand{\N }{\mathbb{N}}
\newcommand{\nfrak}{\mathfrak{n}}
\newcommand{\op }{\mathrm{op}}
\newcommand{\ot }{\otimes }
\newcommand{\pair }[2]{\langle #1,#2 \rangle }

\newcommand{\bigpair }[2]{\Big\langle #1,#2 \Big\rangle }
\newcommand{\ppair }[2]{\langle\!\langle #1,#2 \rangle\!\rangle }
\newcommand{\bigppair }[2]{\Big\langle\!\!\Big\langle #1,#2
                                       \Big\rangle\!\!\Big\rangle }

\newcommand{\pfrak}{\mathfrak{p}}
\newcommand{\pc}{\check{P}}
\newcommand{\ph}{\hat{P}}
\newcommand{\qh}{\hat{Q}}
\newcommand{\qc}{\check{Q}}
\newcommand{\pomE }{\pi_{\Omega,E} }
\newcommand{\pomF }{\pi_{\Omega,F} }
\newcommand{\pord }{\succ }
\newcommand{\pordne }{{\,\scriptstyle \succnsim \,}}

\newcommand{\Rfrak }{\mathfrak{R}}

\newcommand{\rh }{\hat{R}}
\newcommand{\rid }{\mathcal{R}}
\newcommand{\ra }{\acute{R}}
\newcommand{\rg }{\grave{R}}
\newcommand{\rc }{\check{R}}
\newcommand{\rhm }{\hat{R}^-{}}
\newcommand{\ram }{\acute{R}^-{}}
\newcommand{\rgm }{\grave{R}^-{}}
\newcommand{\rcm }{\check{R}^-{}}

\newcommand{\roots }{R }
\newcommand{\Rqps }{\overline{{\roots }^+_S}}

\newcommand{\sqgp}{S_q[G/P_S]}
\newcommand{\sqgpcc}{S_q[G/P_S]_\C^{c=1}}
\newcommand{\sqgpop}{S_q[G/P_S^{{\mathrm{op}}}]}

\newcommand{\TomE }{T_{\Omega,E}}
\newcommand{\TomF }{T_{\Omega,F}}

\newcommand{\U }{U}

\newcommand{\uqb }{U_q(\mathfrak{b})}
\newcommand{\uqg }{U_q(\mathfrak{g})}
\newcommand{\uqks }{U_q(\mathfrak{k}_S)}
\newcommand{\uqls }{U_q(\mathfrak{l}_S)}
\newcommand{\uqgbarplus }{\bar{U}^+_q(\mathfrak{g})}
\newcommand{\Uqbp }{U_q(\mathfrak{b}_+)}
\newcommand{\Uqbm }{U_q(\mathfrak{b}_-)}
\newcommand{\Uqnp }{U_q(\mathfrak{n}_+)}
\newcommand{\Uqnm }{U_q(\mathfrak{n}_-)}
\newcommand{\Ubar }{\overline{\U }}
\newcommand{\wght }{\mathrm{wt}}
\newcommand{\wlat }{P}
\newcommand{\wurz }{\pi }
\newcommand{\vep }{\varepsilon }

\newcommand{\Xlast}{ {}^*\hspace{-.1cm}X }
\newcommand{\Z}{\mathbb{Z}}

\title{De Rham Complex for Quantized Irreducible Flag Manifolds}

\author{Istv\'an Heckenberger and Stefan Kolb}

\date{{\footnotesize\textit{Mathematisches Institut, Universit\"at Leipzig,\\
          Augustusplatz 10,
         04109 Leipzig, Germany}\\Istvan.Heckenberger@math.uni-leipzig.de\qquad
       kolb@itp.uni-leipzig.de
      }\\[\baselineskip] July 31, 2003}

\begin{document}
\maketitle

\begin{abstract}
  It is shown that quantized irreducible flag manifolds possess a canonical
  $q$-analogue of the de Rham complex. Generalizing the well known situation
  for the standard Podle\'s' quantum sphere this analogue is obtained as the
  universal
  differential calculus of a distinguished first order differential
  calculus. The corresponding differential $\dif$ can be written as a
  sum of differentials $\del$ and $\delb$. The universal differential
  calculus corresponding to the first order differential calculi
  $\dif$, $\del$, and $\delb$ are given in terms of generators and relations.
  Relations to well known quantized exterior algebras are established. 
  The dimensions of the homogeneous components are shown to be the same as in
  the classical case. The existence of a volume form is proven. 
\end{abstract}  
{\small

}
\section{Introduction}

The theory of quantum groups provides numerous examples of $q$-deformed
coordinate rings of spaces with group action.
Originally initiated by S.\,L.\, Woronowicz there exists by now a rich theory
of covariant differential calculi over these comodule algebras. 

In A.~Connes' more general concept of noncommutative geometry \cite{b-Connes1}
spectral triples and in particular the Dirac operator are central notions.
One deals with a representation of an algebra $\B$ on some
Hilbert space $\cH$ and with an operator $D:\cH\rightarrow \cH$ such that
the commutators $\dif\,b:=[b,D]$, $b\in \B$, lead to a differential graded
algebra. It has recently been pointed out, that the theory of quantum groups
provides a large class of examples, so called quantized irreducible flag
manifolds, which seem to fit well into Connes' framework of noncommutative
geometry \cite{a-DabSit02p},\cite{a-Kraehmer03p}, \cite{a-SchmWag03p}.
Covariant differential calculi  over quantized irreducible flag manifolds
have been classified in \cite{a-heko03p}. There exists a canonical covariant
first order differential calculus on these spaces, which turned out to
correspond to the Dirac operator constructed in \cite{a-Kraehmer03p}. It is
therefore natural to investigate the corresponding higher order differential
calculi. 

Higher order differential calculi have previously been studied from several
points of view. In \cite{a-SinVa}, \cite{a-SinShklyVa} L.\,L.\,Vaksman and his
coworkers presented a canonical construction of differential forms on quantum
prehomogeneous vector spaces. As these spaces are big cells of the flag
manifolds under consideration here, these differential calculi are closely
related to those investigated in the present paper.

An approach to noncommutative geometry modelled on classical geometry and
compatible with the intrinsic structure of quantum groups and quantum spaces
has been put forward by T.~Brzezi{\'n}ski and S.~Majid in \cite{a-BrzMaj1},
\cite{a-BrzMaj00}.
In this approach the notion of differential calculus is one starting point.
A similar point of view has been adopted by M.~D\Dj ur\dj evi\'c in
\cite{a-Durdevic96}, \cite{a-Durdevic97}.

In the early days of quantum groups there appeared several examples of
differential calculi which in many respects behave as the de Rham complex
over the corresponding commutative algebras \cite{a-Woro3}, \cite{a-PuWo89},
\cite{a-WessZu91}. Yet it became clear that
imposing covariance, i.e.~compatibility with a quantum group action,
one can not expect that such a calculus exists for an arbitrary quantum space.
Nevertheless, there soon existed a well developed theory of covariant
differential calculi on quantum groups while for quantum spaces similar
results have only recently been established \cite{a-Herm02}, \cite{a-HK-QHS}.
Apart from the various examples of differential calculi on quantum groups
(cp.\, references in \cite{b-KS}) and quantum vector spaces as above, only
differential calculi over
Podle\`s' quantum spheres \cite{a-Po92} and Vaksman-Soibelman-spheres
\cite{a-Welk98} have been in detail investigated.

In the present paper differential calculi over quantized irreducible flag
manifolds are studied in detail. Contrary to the situation for quantum groups
and quantum vector spaces for this large class of examples the modules of
differential forms are generally not free over the coordinate algebra. Note,
however, that the general theory implies that being covariant
these modules are projective. The differential calculus constructed here is a
close analogue of the de Rham complex over the corresponding complex
manifold. As in complex geometry the differential $\dif$ can be decomposed
into the sum of differentials $\del$ and $\delb$. The universal differential
calculus corresponding to the first order differential calculi
$\dif$, $\del$, and $\delb$ are given in terms of generators and relations.
The dimensions of the homogeneous components are shown to be the same as in
the classical case. In particular the differential $\dif$ admits a uniquely
determined volume form of degree $2M$, where $M$ is the complex dimension of
the manifold. The fibers of the differential calculi over the classical point
$\vep$ of the quantized flag manifold are shown to be isomorphic to well
known examples of quantized exterior algebras \cite{a-FadResTak1}.

The organization of the paper is as follows. In Chapter
 \ref{preliminaries}
we mainly recall the relevant notions from the theory of quantum groups,
quantum homogeneous spaces, and differential calculus. It is explained in
Section \ref{determining} how the notion of quantum tangent
space introduced in \cite{a-HK-QHS} can be employed to determine the
homogeneous component of degree two of the universal differential calculus
corresponding to a given finite dimensional covariant first order differential
calculus.

Chapter \ref{DCQIFM} is devoted to the construction and
investigation of the desired differential calculus on quantized flag
manifolds. In Section \ref{QFM} the various quantized coordinate rings
associated to flag manifolds are recalled. On the one hand there exist
homogeneous coordinate rings $\sqgp$ and $\sqgpop$. On the other hand
the quantized algebra of functions $\cqgl$ on the quotient $G/L_S$ of the
Lie group $G$ by the Levi factor $L_S$ of the parabolic subgroup
$P_S\subset G$ is considered. It is crucial, as observed in \cite{a-stok02p},
\cite{a-heko03p}, that certain products of generators of $\sqgp$ and $\sqgpop$
generate $\cqgl$. This observation allows the construction of first order
differential calculi $\Gamma_\del$, $\Gamma_\delb$, and $\Gamma_\dif$ over
$\cqgl$ via the construction of first order differential calculi over
$\sqgp$ and $\sqgpop$ in Section \ref{QFFODC}. All first order
differential calculi over $\cqgl$ constructed in this section are also given
in terms of generators and relations, their quantum tangent spaces are
determined, and their dimensions are calculated.

Finally, Section \ref{HODC} is devoted to the corresponding universal
differential calculi $\Gduw$, $\Gdbuw$, and $\Gdifuw$. Again, first the
situation for $\Gamma_\del$ and $\Gamma_\delb$ is analyzed in detail.
Then it is shown that the
differentials $\del$ and $\delb$ can be extended to the universal
differential calculus $\Gdifuw$.
Thus one can reduce statements about $\Gdifuw$ to the
corresponding statements about the the universal
differential calculi $\Gduw$ and $\Gdbuw$.

All algebras considered in this paper are unital $\C$-algebras, likewise
all vector spaces are defined over $\C$.

Throughout this paper several filtrations are defined in the following way.
Let $A$ denote an algebra generated by the elements of a set $Z$ and $\cS$ a
totally ordered abelian semigroup. Then any map
$\deg:Z\rightarrow \cS$ defines a filtration $\Fil$ of the algebra $A$ as
follows.
An element $a\in A$ belongs to $\Fil_n$, $n\in \cS$, if and only if it can be
written
as a polynomial in the elements of $Z$ such that every occurring summand
$a_{1\dots k}z_{1}\dots z_{k}$, $a_{1\dots k}\in\C$, $z_i\in Z$, satisfies
$\sum_{j=1}^k\deg(z_j)\le n$. Instead of $a\in\Fil_n$ by slight abuse of
notation we will also write $\deg(a)=n$.

For any Hopf algebra $H$ the symbols $\kow$, $\vep$, and $\kappa$ will
denote the coproduct, counit, and antipode, respectively. Sweedler notation
for coproducts $\kow a=a_{(1)}\ot a_{(2)}$, $a\in H$, will be used.
If the antipode $\kappa$ is invertible we will frequently identify left and
right $H$-module structures on a vector space $V$ by $vh=\kappa^{-1}(h)v$,
$v\in V$, $h\in H$. The symbol $H^\op$ will denote the corresponding Hopf
algebra with opposite multiplication. 

\section{Preliminaries}\label{preliminaries}
\subsection{Notations}

First, to fix notations some general notions related to Lie algebras are
recalled. Let $\gfrak$ be a finite dimensional complex simple Lie algebra and
$\hfrak\subset \gfrak$ a fixed Cartan subalgebra.
Let $R\subset \hfrak^\ast$ denote the root system
associated with $(\gfrak,\hfrak)$.
Choose an ordered basis $\wurz=\{\alpha_1,\dots,\alpha_r\}$ of simple roots
for $R$ and let $R^+$ (resp.~$R^-$) be the set of positive (resp.~negative)
roots with respect to $\wurz$.
Moreover, let $\gfrak=\nfrak_+\oplus \hfrak\oplus \nfrak_-$ be the
corresponding triangular decomposition. 
Identify $\hfrak$ with its dual via the Killing form. The induced
nondegenerate symmetric bilinear form on $\hfrak^*$
is denoted by $(\cdot,\cdot)$. The root lattice $Q=\Z R$
is contained in the weight lattice $P=\{\lambda\in\hfrak^\ast\,|\,
(\lambda,\alpha_i)/d_i\in\Z\,\forall \alpha_i\in\wurz\}$ where 
$d_i:=(\alpha_i,\alpha_i)/2$. In order to avoid roots of the deformation
parameter $q$ in the following sections we rescale $(\cdot,\cdot)$ such that
$(\cdot,\cdot):P\times P\rightarrow \Z$.

For $\mu,\nu\in \wlat$ we will write $\mu \pord \nu$ if $\mu-\nu $ is a sum
of positive roots and $\mu\pordne \nu$ if $\mu\pord \nu$ and
$\mu\neq \nu$. As usual we define $Q^+:=\{\mu\in Q\,|\, \mu\pord 0\}$.
The height $\hght:Q^+\rightarrow \N_0$
is given by $\hght(\sum_{i=1}^r n_i \alpha_i)=\sum_{i=1}^r n_i$.

The fundamental weights $\omega_i\in\hfrak^\ast$, $i=1,\dots,r$ are
characterized by $(\omega_i,\alpha_j)/d_j=\delta_{ij}$.
Let $P^+$ denote the set of dominant weights, i.~e.~the $\N_0$-span of
$\{\omega_i\,|\,i=1,\dots,r\}$.
Recall that $(a_{ij}):=(2(\alpha_i,\alpha_j)/(\alpha_i,\alpha_i))$ is the
Cartan matrix of $\gfrak$ with respect to $\wurz$.

For $\mu\in P^+$  let $V(\mu)$ denote the uniquely determined finite
dimensional irreducible left $\gfrak$-module with highest weight $\mu$.
More explicitly there exists a nontrivial vector $v_\mu\in V(\mu)$ satisfying
\begin{align}
  E v_\mu=0,\quad H v_\mu=\mu(H)v_\mu
  \qquad \mbox{ for all } H\in\hfrak,\, E\in\nfrak_+.
\end{align}  
For any weight vector $v\in V(\mu)$ let $\wght(v)\in P$ denote the weight of
$v$, i.e.~$Hv=\wght(v)(H)v$. In particular $\wght(v_1)-\wght(v_2)\in Q$ for all weight
vectors $v_1,v_2\in V(\mu)$.

Let $G$ denote the connected simply connected complex Lie group with Lie
algebra $\gfrak$.
For any set $S\subset \wurz$ of simple roots define $R_S^\pm:=\Z S\cap R^\pm$
and $\overline{R_S^\pm}:=R^\pm\setminus R_S^\pm$. 
Let $P_S$ and $P_S^{\mathrm{op}}$ denote the corresponding standard parabolic
subgroups of $G$ with Lie algebra
\begin{align}
  \pfrak_S=\hfrak\oplus\bigoplus_{\alpha\in R^+\cup R^-_S}\gfrak_\alpha,\qquad
  \pfrak^{\mathrm{op}}_S=\hfrak\oplus\bigoplus_{\alpha\in R^-\cup R_S^+}
  \gfrak_\alpha.
\end{align}
Moreover,
\begin{align*}
  \lfrak_S:=\hfrak\oplus\bigoplus_{\alpha\in R^+_S\cup R^-_S}
\gfrak_\alpha
\end{align*}
is the Levi factor of $\pfrak_S$ and $L_S=P_S\cap P_S^\op\subset G$ denotes
the corresponding subgroup. Later on by slight abuse of notation
we will also write $i\in S$ instead of $\alpha_i\in S$.

The generalized flag manifold $G/P_S$ is called irreducible if the
adjoint representation of $\pfrak_S$ on $\gfrak/\pfrak_S$ is irreducible.
Equivalently, $S=\pi\setminus\{\alpha_i\}$ where $\alpha_i$ appears in any
positive root with coefficient at most one. For a complete list of all
irreducible flag manifolds consult e.g.~\cite[p.~27]{b-BastonEastwood}.

\subsection{Quantum Groups}\label{quantumgroups}

\subsubsection{Definition of $\uqg$}
  We keep the notations of the previous section. Let $q\in \C\setminus\{0\}$
  be not a
root of unity. The $q$-deformed universal enveloping algebra $\uqg$ associated
to $\gfrak$
can be defined to be the complex algebra with generators $K_i,K_i^{-1},
E_i,F_i$, $i=1,\dots,r$, and relations 
\begin{align}
\begin{aligned}
&\begin{aligned}
K_iK_i^{-1}&=K_i^{-1}K_i=1,& K_iK_j&=K_jK_i,\\
K_iE_j&=q^{(\alpha_i,\alpha_j)}E_jK_i,& K_iF_j&=q^{-(\alpha_i,\alpha_j)}
F_jK_i,
\end{aligned}&\\
&\begin{aligned}
E_iF_j-F_jE_i&=\delta_{ij}\frac{K_i-K_i^{-1}}{q_i-q_i^{-1}},\\
\sum_{k=0}^{1-a_{ij}}(-1)^k\left(\begin{array}{c}1-a_{ij}\\k \end{array}
        \right)_{q_i}& E_i^{1-a_{ij}-k}E_jE_i^k=0,&i&\neq j,\\
\sum_{k=0}^{1-a_{ij}}(-1)^k\left(\begin{array}{c}1-a_{ij}\\k \end{array}
        \right)_{q_i}&F_i^{1-a_{ij}-k}F_jF_i^k=0,&i&\neq j,
\end{aligned}&
\end{aligned}
\end{align}
where $q_i:=q^{d_i}$ and the $q$-deformed binomial coefficients are defined
by
\begin{align*}
{n\choose k}_q=\frac{[n]_q[n{-}1]_q\dots[n{-}k{+}1]_q}{[1]_q[2]_q\dots[k]_q},
    \qquad [x]_q=\frac{q^x-q^{-x}}{q-q^{-1}}.
\end{align*}
The algebra $\uqg$ obtains a Hopf algebra structure by
\begin{align}
\kopr K_i&= K_i\otimes K_i,&\kopr E_i&= E_i\otimes K_i+1\otimes E_i,&
\kopr F_i&= F_i\otimes 1 + K_i^{-1}\otimes F_i,\nonumber\\
\epsilon(K_i)&=1,&\epsilon(E_i)&=0,& \epsilon(F_i)&=0,\label{hopfstruc}\\
\kappa(K_i)&=K_i^{-1},&\kappa(E_i)&=-E_iK_i^{-1},& \kappa(F_i)&=-K_iF_i.
\nonumber
\end{align}
Let $\Uqnp,\Uqbp,\Uqnm,\Uqbm\subset \uqg$ denote the subalgebras generated
by $\{E_i\,|\,i=1,\dots,r\}$, $\{E_i,K_i,K_i^{-1}\,|\,i=1,\dots,r\}$,
$\{F_i\,|\,i=1,\dots,r\}$, and $\{F_i,K_i,K_i^{-1}\,|\,i=1,\dots,r\}$,
respectively. Moreover, for $\beta\in Q^+$ we will write
$\U_q^\beta(\nfrak_+):=\{x\in \Uqnp\,|\, K_i x
K_i^{-1}=q^{(\beta,\alpha_i)}x\}$ and
$\U_q^\beta(\nfrak_-):=\{x\in \Uqnm\,|\, K_i x
K_i^{-1}=q^{-(\beta,\alpha_i)}x\}$.

\subsubsection{Type 1 Representations}\label{type1}
For $\mu\in P^+$  let $V(\mu)$ denote the uniquely determined finite
dimensional irreducible left $\uqg$-module with highest weight $\mu$.
More explicitly, there exists a highest weight vector
$v_\mu\in V(\mu)\setminus\{0\}$ satisfying
\begin{align}
  E_iv_\mu=0,\quad
  K_iv_\mu=q^{(\mu,\alpha_i)}v_\mu
  \qquad \mbox{ for all } i=1,\dots,r.
\end{align}
A finite dimensional $\uqg$-module $V$ is called \textit{ of type 1} if
$V\cong \bigoplus_i  V(\mu_i)$ is isomorphic to a direct sum of finitely many
$V(\mu_i)$, $\mu_i\in P^+$. The category $\cC$ of $\uqg$-modules of type 1
is a tensor category. By this we mean that $\cC$ contains the trivial
$\uqg$-module $V(0)$ and satisfies
\begin{align}\label{tenscat}
  X,Y\in\mathcal{C} \Rightarrow X\oplus Y,\,X\ot Y,\,X^*\in \mathcal{C}
\end{align}
where $(uf)(x):=f(\kappa(u)x)$ for all $u\in U$, $f\in X^\ast$, $x\in X$.

During subsequent considerations we will meet various natural right
$\uqg$-modules. As indicated at the end of the introduction we will always
endow a right $\uqg$-module $V$ with the left $\uqg$-action defined by
\begin{align*}
  uv:= v\kappa(u)\qquad \forall u\in U, v\in V.
\end{align*}  

\subsubsection{The Braiding}
The category $\cC$ in Section \ref{type1} is a braided category.
Unfortunately the relevant section
in our main reference \cite[8.3.3]{b-KS} lacks notational consistency.
To be able to derive additional properties of the braiding we recall its
construction in some detail.

Recall that the dual pairing $\pair{\cdot}{\cdot}:\Uqbp\times\Uqbm^\op
\rightarrow \C$ of
Hopf algebras \cite[6.3.1]{b-KS} remains non-degenerate when restricted to
$\Uqnp\times\Uqnm^\op$ and satisfies $\pair{a}{b}=0$ for all
$a\in U_q^\mu(\mathfrak{n}_+)$, $b\in U_q^{-\nu}(\mathfrak{n}_-)$,
$\mu\neq\nu$. Let
$C_\beta\in U_q^\beta(\mathfrak{n}_+)\ot U_q^{-\beta}(\mathfrak{n}_-)$
denote the canonical element with respect to $\pair{\cdot}{\cdot}$, i.e.\,
$C_\beta=\sum_i a_i\ot b_i$ where $\{a_i\}$ is a basis of
$U_q^\beta(\mathfrak{n}_+)$ and $\pair{a_i}{b_j}=\delta_{ij}$.
Define 
\begin{align}
\Rfrak&:=\sum_{\beta\in Q^+} (K_\beta\ot 1)(\kappa^{-1}\ot 1)C_\beta&
&\in\uqgbarplus\bar{\ot}\uqgbarplus, \label{Rfrak} \\
\Rfrak^{-1}&:=\sum_{\beta\in Q^+} C_\beta &&\in
\uqgbarplus\bar{\ot}\uqgbarplus\nonumber
\end{align}  
where $\uqgbarplus$ and $\uqgbarplus\bar{\ot}\uqgbarplus$ are defined in
\cite[6.3.3]{b-KS}. There exists an automorphism $\Phi$ of the algebra
$\uqgbarplus\bar{\ot}\uqgbarplus$ such that
\begin{align*}
  \Phi(K_i\ot 1)=K_i\ot 1,&& &\Phi(1\ot K_i)=1\ot K_i,\\
  \Phi(E_i\ot 1)=E_i\ot K_i^{-1},&& &\Phi(1\ot E_i)=K_i^{-1}\ot E_i,\\
  \Phi(F_i\ot 1)=F_i\ot K_i,&& &\Phi(1\ot F_i)=K_i\ot F_i.\\
\end{align*}
One verifies that $\Rfrak$ and $\Phi$ satisfy the properties stated in
\cite[Theorem 8.18]{b-KS}. We suggest first to check the corresponding
properties of $\Rfrak^{-1}$.

For all $V,W \in \cC$ the action of $\Rfrak$ on $V\ot W$ induces a
$\uqg$-module isomorphism
\begin{align}\label{braiding}
  \rh_{V,W}:V\ot W\rightarrow W\ot V,\qquad
  \rh_{V,W}:=\tau\circ B_{V,W}(\Rfrak(v\ot w))
\end{align}
where $\tau$ denotes the twist $\tau(v\ot w)=w\ot v$, and
$B_{V,W}(v\ot w)=q^{(\mu,\nu)}v\ot w$ for weight vectors $v\in V$ and
$w\in W$ of weight $\mu$ and $\nu$, respectively. The family $(\rh_{VW})$
defines a braiding in $\cC$. To simplify notation we will also write
$\rh_{\mu,\nu}:=\rh_{V(\mu),V(\nu)}$ if $\mu,\nu\in P^+$.

\begin{remark}\upshape{
  An explicit formula for $\Rfrak$ is given in
  \cite[8.3.3, above Thm.~18]{b-KS}.
  Note that the braiding $\rh$ is uniquely determined if one demands that
  (\ref{braiding}) is a $\uqg$-module homomorphism satisfying
  \begin{align}\label{braidingFO}
    \rh_{V,W}(v\ot w)=q^{(\wght(v),\wght(w))}w\ot v +\sum_i w_i\ot v_i
  \end{align}
  where $\wght(w) \pord \wght(w_i)$ and $\wght(v_i) \pord \wght(v)$.
  Indeed, if $v_{\max}\in V$ is a highest weight vector then
  $\rh_{V,W}(v_{\max}\ot w)$ is uniquely determined by (\ref{braidingFO})
  for any $w\in W$. The property of being a $\uqg$-module homomorphism
  then fixes $\rh_{V,W}$ on all of $V\ot W$.

  For this reason (\ref{Rfrak}) should coincide with the
  explicit expression in \cite[8.3.3]{b-KS}. It is straightforward to check
  that coefficients of the terms $E_i\ot F_i$ of the two expressions are
  identical. Coefficients of higher order terms of the explicit expression
  of $\Rfrak$ will not be used in this paper.}
\end{remark}  

\subsubsection{Restriction of the Braiding to $\uqls$}
Let $S\subset \pi$ and let $\kfrak_S:=[\lfrak_S,\lfrak_S]\subset\gfrak$
denote the semisimple part of
$\lfrak_S\subset\pfrak_S\subset \gfrak$. Define $\uqks$ and $\uqls$ to be the
Hopf subalgebras of $\uqg$ generated by the sets 
$\{K_j,K_j^{-1},E_j,F_j\,|\, j\in S\}$ and
$\{K_i,K_i^{-1},E_j,F_j\,|\, j\in S,\, i=1,\dots,r\}$, respectively. 

As above the tensor category
$\cC^\kfrak$ of type 1 representations of $\uqks$ is braided with braiding
$\rh^\kfrak_{V,W}$. Moreover, let $(\cdot,\cdot)_\kfrak$ denote the uniquely
determined bilinear form on the weight lattice corresponding to $\kfrak_S$
such that $(\alpha,\beta)_\kfrak=(\alpha,\beta)$ for all simple roots
$\alpha,\beta\in S$. 

The following Lemma will be used only in the proof of Proposition
\ref{hodProp} at the very end of this paper.
\begin{lemma}\label{restrictedR}
  Let $V=V(\nu)$ and $W=V(\mu)$ be irreducible $\uqg$-modules and let
  $V'=V(\nu')\subset V$ and $W'=V(\mu')\subset W$ be irreducible
  $\uqls$-submodules. Let $p_V:V\rightarrow V'$ and $p_W:W\rightarrow W'$
  denote surjective $\uqls$-module homomorphisms satisfying $p_V^2=p_V$ and
  $p_W^2=p_W$, respectively. Then
  \begin{align*}
    (p_W\ot p_V)\rh_{V,W}|_{V'\ot W'}=
    q^{(\nu-\gamma_1,\mu-\gamma_2)-(\nu',\mu')_\kfrak} \rh^\kfrak_{V',W'}
  \end{align*}
  where the highest weight vectors of
  $V'$ and $W'$ have weight
  $\nu-\gamma_1$ in $V$ and $\mu-\gamma_2$ in $W$, $\gamma_1,\gamma_2\in Q^+$,
  respectively. 
\end{lemma}
\begin{proof}
  In analogy to (\ref{Rfrak}) one has an element
  \begin{align*}
    \Rfrak^\kfrak:=\sum_{\beta\in \Z S\cap Q^+}
  (K_\beta\ot 1)(\kappa^{-1}\ot 1)C_\beta \in
  \bar{U}^+_q(\kfrak_S)\bar{\ot}\bar{U}^+_q(\kfrak_S)
  \end{align*}
  and linear maps
  \begin{align*}
    B^\kfrak_{V',W'}:V'\ot W'\rightarrow V'\ot W',\quad v\ot w\mapsto
    q^{(\mu,\nu)_\kfrak}v\ot w
  \end{align*}
  where $V',W'\in \cC^\kfrak$ and $v$ and $w$ are weight vectors of
  weight $\mu$ and $\nu$, respectively.
  Then for all $v'\in V'$, $w'\in W'$ one obtains
  \begin{align*}
    (p_W\ot p_V)\rh_{V,W}(v'\ot w')&\overset{\small (\ref{braiding})}{=}
    \tau\circ B_{V,W} (p_V\ot p_W)
    (\Rfrak(v'\ot w'))\\
    &\overset{\small (\ref{Rfrak})}{=}
    \tau\circ B_{V,W}(\Rfrak^\kfrak(v'\ot w'))\\
    &=B_{W,V}\circ(B^\kfrak_{W',V'})^{-1}\circ\tau\circ B^\kfrak_{V',W'}
      (\Rfrak^\kfrak(v'\ot w'))\\
    &=B_{W,V}\circ(B^\kfrak_{W',V'})^{-1}\circ \rh^\kfrak_{V',W'}(v'\ot w').
  \end{align*}
  Note that by definition of $\gamma_1$ and $\gamma_2$ one has
  $(\nu{-}\gamma_1,\alpha_i)=(\nu',\alpha_i)_\kfrak$ and
  $(\mu{-}\gamma_1,\alpha_i)=(\mu',\alpha_i)_\kfrak$ for all $\alpha_i\in S$.
  Now the claim of the lemma follows from
  \begin{align*}
    (\nu{-}\gamma_1{-}\beta_1,&\mu{-}\gamma_2{-}\beta_2)-
    (\nu'{-}\beta_1,\mu'{-}\beta_2)_\kfrak\\
    =&(\nu{-}\gamma_1,\mu{-}\gamma_2)-
    (\beta_1,\mu{-}\gamma_2)-(\nu{-}\gamma_1,\beta_2)+(\beta_1,\beta_2)\\
    &-(\nu',\mu')_\kfrak+
    (\beta_1,\mu')_\kfrak+(\nu',\beta_2)_\kfrak-(\beta_1,\beta_2)_\kfrak\\
    =&(\nu{-}\gamma_1,\mu{-}\gamma_2)-(\nu',\mu')_\kfrak
  \end{align*}
    for all $\beta_1,\beta_2\in \Z S\cap Q^+$.
\end{proof}  

\subsubsection{$R$-Matrices}\label{subsRMatrices}
To write coordinate algebras of quantized flag manifolds in terms of
generators and relations it will be helpful to introduce
additional notations for certain special cases of $\rh$.
For $\lambda=\sum_{i\notin S}\omega_i$
set $N:=\dim V(\lambda)$ and abbreviate $I:=\{1,\dots,N\}$. Choose a basis
$\{v_i\,|\,i \in I\}$ of weight vectors of $V(\lambda)$ and
let $\{f_i\,|\,i\in I\}$ be the corresponding dual basis.
Define matrices $\rh$, $\rc$, $\ra^-$ and $\rg^-$ by
\begin{align*}
  \rh_{\lambda,\lambda}(v_i{\ot} v_j)&=:
                     \sum_{k,l\in I}\rh^{kl}_{ij} v_k{\ot} v_l,&
  \rh_{-w_0\lambda,-w_0\lambda}(f_i{\ot} f_j)&=:
                     \sum_{k,l\in I}\rc^{kl}_{ij} f_k{\ot} f_l,\\
  \rh_{-w_0\lambda,\lambda}(f_i{\ot} v_j)&=:
                     \sum_{k,l\in I}\ram {}^{kl}_{ij} v_k{\ot} f_l,&
  \rh_{\lambda,-w_0\lambda}(v_i{\ot} f_j)&=:\sum_{k,l\in I}
       \rgm {}^{kl}_{ij} f_k{\ot} v_l.
\end{align*}  
Alternatively
\begin{align*}
  (f_i{\ot} f_j){\circ} \rh_{\lambda,\lambda}&=
                       \sum_{k,l\in I}\rh_{kl}^{ij} f_k{\ot} f_l,&
  (v_i{\ot} v_j){\circ} \rh_{-w_0\lambda,-w_0\lambda}&=\sum_{k,l\in I}
                       \rc_{kl}^{ij} v_k{\ot} v_l,\\
  (f_i{\ot} v_j){\circ} \rh_{-w_0\lambda,\lambda}&=
                       \sum_{k,l\in I}\ram {}_{kl}^{ij} v_k{\ot} f_l,&
  (v_i{\ot} f_j){\circ} \rh_{\lambda,-w_0\lambda}&=
                       \sum_{k,l\in I}\rgm {}_{kl}^{ij} f_k{\ot} v_l,\\
\end{align*}  
where the elements of $V(\lambda)$ are considered as functionals on
$V(\lambda)^\ast$. Let $\rh^-$, $\rc^-$, $\rg$ and $\ra$ denote the
inverse of the matrix $\rh$, $\rc$, $\ra^-$ and $\rg^-$, respectively.

By (\ref{Rfrak}) the matrix $\rh $ has the property that
$\rh ^{ij}_{kl}\not=0$
implies that $i=l,j=k$ or both $\wght (v_j)\pordne\wght (v_k)$ and
$\wght (v_l)\pordne\wght (v_i)$. Therefore we associate to $\rh $ the symbol
$<$ which denotes the positions of the larger weights.
Similar properties are fulfilled for the other types of $R$-matrices.
For example, the relation $\ram {}^{ij}_{kl}\neq 0$ implies that
$i=l,j=k$ or both
$\wght (v_k)\pordne\wght (v_j)$ and $\wght (v_l)\pordne\wght (v_i)$.
We collect these properties in the following table.

\begin{gather}\label{rtabelle}
\begin{array}{|c|c|c|c|c|c|c|c|}
\hline
\rh \rule{0pt}{3ex} & \rhm & \ra & \ram & \rc & \rcm & \rg & \rgm \\
\hline
< & > & \vee & \wedge & > & < & \wedge & \vee \\
\hline
\end{array}
\end{gather}

\subsubsection{The $q$-Deformed Coordinate Ring $\cqg$}

The $q$-deformed coordinate ring $\cqg$ is defined to be the subspace
of the linear dual $\uqg^\ast$ spanned by the matrix coefficients of the
finite dimensional irreducible representations $V(\mu)$, $\mu\in
P^+$.
For $v\in V(\mu)$, $f\in V(\mu)^\ast$ the matrix coefficient
$c^\mu_{f,v}\in \uqg^\ast$ is defined by
\begin{align*}
  c^\mu_{f,v}(X)=f(Xv).
\end{align*}
The linear span of matrix coefficients of $V(\mu)$
\begin{align}
  C^{V(\mu)}=\Lin_\C\{c^\mu_{f,v}\,|\,v\in V(\mu), f\in V(\mu)^\ast\}
\end{align}
obtains a $\uqg$-bimodule structure by
\begin{align}
  (Yc^\mu_{f,v}Z)(X)=f(ZXYv)=c^\mu_{fZ,Yv}(X).
\end{align}
Here $V(\mu)^\ast$ is considered as a right $\uqg$-module. Note that by
construction
\begin{align}
  \cqg\cong\bigoplus_{\mu\in P^+} C^{V(\mu)}
\end{align}
is a Hopf algebra and the pairing
\begin{align}\label{qpair}
  \cqg\ot \uqg\rightarrow \C
\end{align}
is nondegenerate.

\subsubsection{Quantum Homogeneous Spaces}\label{qhs}
We recall the class of quantum homogeneous spaces considered in
\cite{a-MullSch}.
Let $\U $ denote a Hopf algebra over $\C$ with bijective antipode $\kappa$ and
$\K \subset \U $ a right coideal subalgebra with right coaction
$\Kkow :\K \rightarrow \K \ot \U $. Consider a tensor category
$\mathcal{C}$ of finite dimensional left $\U $-modules.
By this we mean that $\mathcal{C}$ is a class of finite
dimensional left $\U$-modules containing the trivial $\U$-module via
$\vep $ and  satisfying (\ref{tenscat}).

Let $\cA :=\U ^\circ_{\mathcal{C}}$ denote the dual Hopf algebra generated by
the matrix coefficients of all $\U $-modules in $\mathcal{C}$.
Assume that $\cA$ separates the elements of $\U$.
Assume further that the antipode of $\cA$ is bijective. Note that this is
equivalent to
\[
  X\in\mathcal{C}\Rightarrow {}^*\hspace{-.1cm}X \in \mathcal{C}
\]
where $\Xlast=X^*$ as a vector space and $(uf)(x):=f(\kappa^{-1}(u)x)$ for all
$u\in U$, $f\in \Xlast$, $x\in X$.

Define a left coideal subalgebra $\B \subset \cA $ by
\begin{equation}\label{Bdef}
  \B :=\{b\in \cA \,|\, b_{(1)} \,b_{(2)}(k)=\vep(k)b \quad
  \mbox{for all $k\in \K$}\}.
\end{equation}
Assume $\K $ to be $\mathcal{C}$-semisimple, i.e. the restriction of any
$\U$-module in $\mathcal{C}$ to the subalgebra $K\subset \U$ is isomorphic
to the direct sum of irreducible $K$-modules. In full analogy to
\cite[Thm.~2.2 (2)]{a-MullSch} this implies that $\cA $ is a faithfully
flat $\B $-module.

\subsubsection{Categorial Equivalence}

Assume $\B\hookrightarrow \cA$ to be a left coideal subalgebra
of a Hopf algebra $\cA$ with bijective antipode and define
$\B^+:=\{b\in \B\,|\, \vep(b)=0\}$. Then $\Aqr:=\cA/\B^+\cA$ and
$\Aql:=\cA/\cA\B^+$ are right and left $\cA$-module coalgebras, respectively.
Let $\ABM$ and $\AqM$ denote the category of left $\cA$-covariant
left $\B$-modules and of left $\Aqr$-comodules, respectively.
Recall that for any coalgebra $C$ the cotensor product of a left
$C$-comodule $P$ and a right $C$-comodule $Q$ is defined by
\begin{align*}
  P\square_C Q :=\left\{\sum_ip_i{\ot} q_i\in P\ot Q\,\bigg|\,\sum_i
   p_{i(0)}{\ot} p_{i(1)}{\ot} q_i=\sum_i p_i{\ot} q_{i(-1)}{\ot}q_{i(0)}
                \right\}.
\end{align*}  
There exist functors
\begin{align*}
  \Phi:\ABM&\rightarrow \AqM,& \Phi(\Gamma)&=\Gamma/\B^+\Gamma,\\
  \Psi:\AqM&\rightarrow \ABM,& \Psi(V)&=\cA\square_{\Aqr} V.
\end{align*}
Here for any $\Gamma\in \ABM$ the left $\Aqr$-comodule structure
on $\Gamma/\B^+\Gamma$ is induced by the left $\cA$-comodule structure 
of $\Gamma$. Moreover, the left $\B$-module and the left
$\cA$-comodule structures of $\cA\square_{\Aqr} V$ are defined on the
first tensor factor.

\begin{theorem}{\em \cite[Theorem 1]{a-Tak79}}\label{catequiv}
  With the notions as above suppose that $\cA$ is a faithfully flat
  right $\B$-module. Then $\Phi$ and $\Psi$ are mutually inverse
  equivalences of categories.
\end{theorem}  
\begin{remark}\label{catequivremarks}\upshape{
(i) The functor $\Phi$ is equivalent to $\Phi':\ABM\rightarrow \AqM$ defined
by
\begin{align*}
  \Phi'(\Gamma)&:={}^{\co \cA}(\cA\ot_\B\Gamma)\\
               &:=\left\{\sum_i a_i\ot \rho_i\,\bigg |\, \sum_i a_{i(1)}
               \rho_{i(-1)}\ot a_{i(2)}\ot\rho_{i(0)}=1\ot \sum_i
               a_i\ot\rho_i\right \}
\end{align*}
where the left $\Aqr$-comodule structure on ${}^{\co \cA}(\cA\ot_\B\Gamma)$
is given by
\begin{align*}
  \kow_L\left(\sum_i a_i\ot \rho_i\right)=
  \sum_i \kappa^{-1}(a_{i(2)})\ot (a_{i(1)}\ot \rho_i).
\end{align*}
The isomorphism ${}^{\co \cA}(\cA\ot_\B\Gamma)\rightarrow \Gamma/\B^+\Gamma$ is
defined by $\sum_i a_i\ot \rho_i\mapsto \sum_i\vep(a_i)\overline{\rho_i}$.
To verify injectivity note that
${}^{\co \cA}(\cA\ot_\B\Gamma)=
  \{\kappa(\rho_{(-1)})\ot\rho_{(0)}\,|\,\rho\in \Gamma\}.$

(ii) The functor $\Psi$ is equivalent to $\Psi':\AqM\rightarrow \ABM$ defined
   by
   \begin{align*}
     &\Psi'(V):=(\cA\ot V)^{\co \Aql}\\
          &:=\left\{ \sum_i a_i{\ot} v_i\in \cA{\ot} V\,\bigg|\,
          \sum_i a_{i(1)}{\ot} v_{i(0)}
          {\ot} a_{i(2)}\kappa(v_{i(-1)})=\sum_i a_i {\ot} v_i{\ot} 1\right\}.
   \end{align*}
Here the left $\B$-module and left $\cA$-comodule structure on $\Psi'(V)$ is
given by
\begin{align*}
  b\left(\sum_i a_i{\ot} v_i\right)=\sum_i (ba_i)\ot v_i,\quad
  \kow_L\left(\sum_i a_i{\ot} v_i\right)=\sum_i a_{i(1)}\ot (a_{i(2)}\ot v_i)
\end{align*}
for all $b\in \B$ and $\sum_i a_i\ot v_i\in (\cA\ot V)^{\co \Aql}$.

(iii) In the situation of Theorem \ref{catequiv} the coalgebra $\Aqr$ is
 cosemisimple. Therefore
  any $\Gamma\in \ABM$ is a projective left $\B$-module.
}\end{remark}    
  
\subsection{Differential Calculus}
\subsubsection{First Order Differential Calculus}
 For the convenience of the reader the notion of differential calculus from
\cite{a-Woro2} is recalled. 
A \textit{first order differential calculus} (FODC)
over an algebra $\B$ is a $\B$-bimodule $\Gamma$ together with a
$\C$-linear map
\begin{equation*}
  \dif:\B\rightarrow\Gamma
\end{equation*}
such that $\Gamma=\Lin_\C\{a\,\dif b\,c\,|\,a,b,c\in\B\}$ and $\dif$
satisfies the Leibniz rule
\begin{align*}
  \dif(ab)&=a\,\dif b + \dif a\,b.
\end{align*}    
Let in addition $\cA$ denote a Hopf algebra and
$\kow_\B:\B\rightarrow \cA\otimes\B$ a left $\cA$-comodule algebra
structure on $\B$.
If $\Gamma$ possesses the structure of a left $\cA$-comodule
\begin{equation*}
  \kow_\Gamma:\Gamma\rightarrow\cA\ot \Gamma
\end{equation*}
such that
\begin{equation*}
\kow_\Gamma(a\dif b\,c)=(\kow_\B a)((\id\otimes\dif)\kow_\B b)
                          (\kow_\B c)
\end{equation*}
then $\Gamma$ is called \textit{left covariant}.

For a family of left covariant FODC
$(\Gamma_i,\dif_i)_{i=1,\dots,k}$ define $\dif=\bigoplus_i \dif_i:\B
\rightarrow \bigoplus_i\Gamma_i$. Then
$\Gamma=\B\dif \B\subset \oplus_i \Gamma_i$ is a covariant FODC with
differential $\dif$
which is called the \textit{sum} of the calculi $\Gamma_1,\dots,\Gamma_k$.

If $\cD\subset \B$ is a subalgebra and $(\Gamma,\dif)$ is a FODC over $\B$ then
$(\Gamma|_\cD,\dif|_\cD)$ defined by
\begin{align*}
  \Gamma|_\cD:=\{a\dif b\,|\, a,b\in \cD\},\quad \dif|_\cD(a):=\dif a \quad
  \forall a\in \cD,
\end{align*}  
is a FODC over $\cD$ called the \textit{FODC over $\cD$ induced by $\Gamma$}.

\subsubsection{Higher Order Differential Calculus}
A \textit{differential calculus} (DC) over $\B$ is a differential graded
algebra $(\Gamma^\wedge=\oplus_{i\in \N_0}\Gamma^{\wedge i},\dif)$ such that
$\Gamma^{\wedge 0}=\B$ and $\Gamma^\wedge$ is generated by $\B$ and $\dif
\B$. The product of a DC will usually be denoted by $\wedge$.
Assume that $\B$ is a left $\cA$-comodule algebra over a Hopf algebra $\cA$.
Then a DC $(\Gamma^\wedge,\dif)$ over $\B$ is called (left)
\textit{covariant} if $\Gamma^\wedge$ has the structure of a (left)
$\cA$-comodule algebra such that
$(\Gamma^{\wedge 1},\dif|_\B)$ is a (left) covariant FODC.

The universal DC of a FODC $(\Gamma,\dif_\Gamma)$ over $\B$ is the uniquely
determined DC $(\Gamma_u^\wedge,\dif_u)$ over $\B$ with
$\Gamma_u^{\wedge 1}=\Gamma$ and $\dif_u|_\B=\dif_\Gamma$ satisfying the
following defining property. For any DC $(\Gamma^\wedge,\dif)$ over $\B$ with
$\Gamma^{\wedge 1}=\Gamma$ and $\dif|_\B=\dif_\Gamma$ there exists a map
$\phi:\Gamma^\wedge_u\rightarrow \Gamma^\wedge$ of differential graded
algebras such that $\phi|_{\B\oplus\Gamma}=\id$.  To construct
$(\Gamma_u^\wedge,\dif_u)$ consider the tensor algebra
$\Gamma^\ot=\bigoplus_{k=0}^{\infty}\Gamma^{\ot k}$ of the $\B$-bimodule
$\Gamma$.
Then $\Gamma_u^{\wedge}$ is the quotient of $\Gamma^\ot$ by the ideal
generated by $\{\sum_i\dif a_i\ot \dif b_i\,|\,\sum_i a_i\dif b_i=0\}$
and the differential is defined by $\dif_u(a_0\dif a_1\wedge \dots \wedge \dif a_n)=
\dif a_0\wedge \dif a_1\wedge\dots \wedge \dif a_n$.

\subsubsection{Right Ideals and Quantum Tangent Spaces}
From now on we assume $K\subset U$ to be a right coideal subalgebra and
$\B\subset \cA=\U_\cC^\circ$ to be the corresponding quantum homogeneous
space as in Subsection \ref{qhs}.

In this situation left covariant first order differential calculi
over $\B$ are in one-to-one correspondence to right ideals
$\rid\subset\B^+$ satisfying $\kow_\B\rid\subset \cA\ot \rid + \B^+\cA\ot\B$
\cite{a-Herm02}. The right ideal corresponding to a covariant FODC $\Gamma$
is given by
\begin{align}\label{rid}
  \rid=\Big\{\sum_i \vep(a_i)b_i^+\,\Big|\,\sum_i a_i\,\dif b_i=0\Big\}
  \subset \B^+
\end{align}  
where $b^+=b-\vep(b)$ for all $b\in\B$.
Conversely, to construct the FODC $\Gamma$ corresponding to $\rid$ consider the
$\B$-bimodule structure on $\tilde{\Gamma}:=\cA\ot (\B^+/\rid)$ given by
\begin{equation}
  c(a\ot \overline{b})c'= cac'_{(-1)}\ot\overline{bc'_{(0)}},\qquad c,c'\in
  \B,\,b\in \B^+,\, a\in \cA
\end{equation}  
and the differential $\dif:\B\rightarrow \tilde{\Gamma}$,
$\dif b={b_{(-1)}}\ot \overline{b_{(0)}^+}$.
Then one obtains $\Gamma$ by 
$\Gamma=\Lin_\C \{ b_1\,\dif b_2\,|\,b_1,b_2\in\B\}$.
This implies in particular 
\begin{align}\label{adb=0}
\sum_i a_i\dif b_i=0\,\Leftrightarrow \, \sum_i a_i b_{i(-1)}\ot b_{i(0)}^+
\in \cA\ot \rid.
\end{align}
 
To a FODC $\Gamma$ with corresponding right ideal $\rid$ one associates the
vector space
\begin{align*}
T^\vep_\Gamma=\{f\in \B^\ast\,|\, f(x)=0 \textrm{ for all }x\in\rid\}
\end{align*}
and the so called \textit{quantum tangent space}
\begin{align*}
  T_\Gamma=(T_\Gamma^\vep)^+=\{f\in T_\Gamma^\vep\,|\, f(1)=0\}.
\end{align*}  
The dimension of a first order differential calculus is defined by
\begin{align*}
  \dim \Gamma=\dim_\C \Gamma/\B^+\Gamma=\dim_\C\B^+/\rid.
\end{align*} 
Let $\B^\circ$ denote the dual coalgebra of $\B$.
\begin{proposition}\label{corresp}{\em \cite[Cor.~5]{a-HK-QHS}}
Let $\B\subset \cA$ be as in Subsection \ref{qhs}. Then there is a canonical
one-to-one correspondence
between $n$-dimensional left covariant FODC over $\B $ and $(n+1)$-dimensional
subspaces $T^\vep\subset \B ^\circ $ such that 
\begin{align}
\vep \in T^\vep,\quad \kow T^\vep \subset  T^\vep \ot \B^\circ,\quad K T^\vep
 \subset T^\vep.
\end{align}
\end{proposition}
A covariant FODC $\Gamma\neq \{0\}$ over $\B$ is called
\textit{irreducible} if it does not possess any nontrivial quotient
(by a left covariant $\B$-bimodule). Note that for finite dimensional calculi
this property is equivalent to the property that $T^\vep_\Gamma$ does not
possess any 
left $K$-invariant right $\B^\circ$-subcomodule $\tilde{T}$ such that
$\C\cdot\vep\varsubsetneq\tilde{T} \varsubsetneq T^\vep_\Gamma$.

Let $\Gamma$ be a sum of finite dimensional covariant FODC
$\Gamma_i$, $i=1,\dots,N$, over $\B$
with corresponding right ideals $\rid_i$.
Then the right ideal corresponding to $\Gamma$ is given by
$\rid_\Gamma=\cap_i \rid_{\Gamma_i}$ and therefore the relation
$T_\Gamma=T_{\Gamma_1}+\dots+T_{\Gamma_k}$ of quantum tangent spaces holds.
A sum of covariant differential calculi is called a \textit{direct sum} if
$\Gamma=\oplus_i \Gamma_i$ is a direct sum of bimodules.
This condition is equivalent to $T_\Gamma=\oplus_i T_{\Gamma_i}$.

\subsubsection{Induced Covariant FODC}
Using quantum tangent spaces it is possible to identify induced covariant
FODC.
\begin{proposition}{\em \cite[Cor.~9]{a-HK-QHS}}\label{induced}
Let $\B\subset \cA$ be as in Subsection \ref{qhs} and let $\Gamma$ be a
finite dimensional left-covariant FODC over $\cA$ with quantum tangent space
$T$. Then $\Gamma|_{\B}$ is finite dimensional if and only if $K T|_{\B}$ is
finite dimensional. In this case the quantum tangent space of
$\Gamma|_{\B}$ coincides with $K T|_{\B}$.
\end{proposition}

\begin{lemma}\label{equivalences}
  Let $\B\subset \cA$ be as in Subsection \ref{qhs} and
  let $\Omega$ denote a finite dimensional covariant FODC over $\cA$ with
  corresponding right ideal $\rid_\Omega$ and quantum tangent space
  $T_\Omega$. Assume that the induced FODC $\Gamma$ over $\B$ is finite
  dimensional with right ideal $\rid$ and quantum tangent space $T$. Then
  the following properties are equivalent:
  \begin{itemize}
    \item[(i)] $T=T_\Omega|_\B$.
    \item[(ii)] $\rid=\rid_\Omega\cap \B$.
    \item[(iii)] The canonical map $\cA\ot_\B\Gamma \rightarrow \Omega$,
      $a\ot \gamma\mapsto a\gamma$ is injective.
  \end{itemize}  
\end{lemma}  
\begin{proof}
  The canonical map $\cA\ot_\B\Gamma \rightarrow \Omega$ of left covariant
  left $\cA$-modules is injective if and only if the restriction
  \begin{align}\label{invinclusion}
  {}^{\co \cA}(\cA\ot_\B\Gamma) \rightarrow {}^{\co \cA}\Omega
  \end{align}
  to the space of left
  coinvariants is injective.
  Recall from \cite[Thm.~5.1]{a-Woro2}, \cite[Lemma 6]{a-HK-QHS} that there
  exist pairings
  \begin{align}
    \pair{\cdot}{\cdot}_\Omega&:\Omega\times T_\Omega\rightarrow \C,&
    \pair{a\dif b}{t}_\Omega &=\vep(a)t(b),\label{OmTpairing}\\
    \pair{\cdot}{\cdot}  &:\Gamma \times T \rightarrow \C,&
    \pair{a\dif b}{t} &=\vep(a)t(b),\label{GamTpairing} 
  \end{align}  
  which induce nondegenerate pairings
  \begin{align}\label{invTpairings}
        {}^{\co \cA}\Omega\times T_\Omega\rightarrow \C,\qquad
        {}^{\co \cA}(\cA\ot_\B \Gamma)\times T \rightarrow \C.  
  \end{align}
  Now, (\ref{invinclusion}) is injective if and only if $T_\Omega$ separates
  ${}^{\co \cA}(\cA\ot_\B\Gamma)$. In view of the nondegeneracy of the second
  pairing in (\ref{invTpairings}) and Proposition
  \ref{induced} the latter is equivalent to $T=T_\Omega|_\B$.
  Therefore (i) is equivalent to (iii). The equivalence
  between (i) and (ii) holds by duality.
\end{proof}

\subsubsection{Determining $\Gamma^{\wedge 2}_u$}\label{determining}
Quantum tangent spaces can also be employed to obtain information about
  higher order differential calculi.
  Let $\Omega$, $\Gamma$, $T_\Omega$, and $T$ be as in
  Lemma \ref{equivalences}. In analogy to
  (\ref{OmTpairing}), (\ref{GamTpairing})  there exists a pairing
  \begin{align}
    &\ppair{\cdot}{\cdot}:
    \cA\ot_\B\Gamma\ot_\B\Gamma \times T_\Omega\ot T\rightarrow \C,\nonumber\\
    &\ppair{ a\ot \rho\ot \zeta}{s\ot t}:=\vep(a)\pair{\rho}{s_{(0)}^+} _\Omega
    \,s_{(1)}(\zeta_{(-1)})\pair{\zeta_{(0)}}{t}.
    \label{ppairingdef}
  \end{align}
  In particular
  \begin{align}\label{adbdcppair}
     \ppair{ a\ot \dif b\ot \dif c}{s\ot t}=\vep(a)s(b^+c_{(-1)})t(c_{(0)}).
  \end{align}  
  To verify that $\ppair{\cdot}{\cdot}$ is well defined note that
  \begin{align*}
     \pair{a\dif b\,c}{s}_\Omega=\pair{a\dif(bc)-ab\dif c}{s}_\Omega
     =\vep(a)s(b^+c)=\pair{a\dif b}{s_{(0)}^+}_\Omega s_{(1)}(c)
  \end{align*}
  and therefore
  \begin{align*}
     \pair{\rho c}{s_{(0)}^+}_\Omega\,s_{(1)}(\zeta_{(-1)})
     \pair{\zeta_{(0)}}{t}&=
     \pair{\rho }{s_{(0)}^+}_\Omega\,s_{(1)}(c)s_{(2)}(\zeta_{(-1)})
     \pair{\zeta_{(0)}}{t}\\
     &= \pair{\rho }{s_{(0)}^+}_\Omega\,s_{(1)}(c_{(-1)}\zeta_{(-1)})
     \pair{c_{(0)}\zeta_{(0)}}{t}.
  \end{align*}  
  \begin{lemma}\label{ppair}
     Let $\B\subset \cA$ be as in Subsection \ref{qhs} and let
     $\Omega$, $\Gamma$, $T_\Omega$, and $T$ be as in
     Lemma \ref{equivalences}. Assume $T=T_\Omega|_\B$. Then the pairing
     $\ppair{\cdot}{\cdot}$ induces a nondegenerate pairing
     \begin{align}\label{nondegppairing}
       {}^{\co \cA}(\cA\ot_\B\Gamma\ot_\B\Gamma) \times (T_\Omega\ot T)/T_0
       \rightarrow \C, 
     \end{align}
     where $T_0=\{\sum_i s_i\ot t_i \in T_\Omega\ot T\,|\,
                  \sum_i s_{i(0)}^+|_\B\ot s_{i(1)}t_i=0
                  \mbox{ in }\B^\circ\ot\B^\circ \}$.
  \end{lemma}  
  \begin{proof}
     Note first that
     \begin{align*}
        \ppair{a\ot\dif b\ot \dif c}{s\ot t}
        =\ppair{\kappa(a_{(1)}b_{(-1)}c_{(-1)})a_{(2)}\ot
               \dif b_{(0)}\ot \dif c_{(0)}}{s\ot t}.
     \end{align*}  
     By definition (\ref{ppairingdef}) the pairing (\ref{nondegppairing}) is
     well defined and by (\ref{adbdcppair}) the elements of $(T_\Omega\ot
     T)/T_0$ are separated by ${}^{\co \cA}(\cA\ot_\B\Gamma\ot_\B\Gamma)$.

     Conversely, recall that $\Gamma$ is a projective left $\B$-module by
     Remark \ref{catequivremarks}(iii).
     By Lemma \ref{equivalences} one obtains a canonical inclusion
     \begin{align*}
        \cA\ot_\B\Gamma\ot_\B\Gamma\subset \Omega\ot_\B\Gamma \cong
        \Omega\ot_\cA\cA\ot_\B\Gamma\subset \Omega\ot_\cA\Omega.
     \end{align*}
     Via this inclusion one identifies
     $\kappa(a_{(-1)}b_{(-1)})\ot\dif a_{(0)}\ot
               {}^{\co \cA}\dif b_{(0)}\in(\cA\ot_\B\Gamma\ot_\B\Gamma)$ with
     $\kappa(b_{(-1)})\omega_L(a)\ot \dif b_{(0)}\in \Omega\ot_\B\Gamma$
     and with
     $\omega_L(a^+b_{(-1)})\ot \omega_L(b_{(0)})\in \Omega\ot_\cA\Omega$
     where $\omega_L(a)=\kappa(a_{(1)})\dif a_{(2)}$ for all $a\in \cA$.

     Recall \cite[p.164]{a-Woro2} that the pairing
     \begin{align*}
        \ppair{\cdot}{\cdot}:{}^{\co \cA}(\Omega\ot\Omega)\times
         (T_\Omega\ot_\C T_\Omega)&\rightarrow \C,\\
         \ppair{\omega_L(a)\ot\omega_L(b)}{s\ot t}&=s(a)t(b),\quad a,b\in \cA
     \end{align*}  
     is nondegenerate and compatible with (\ref{nondegppairing}).
     Therefore
     \begin{align*}
       \bigppair{\sum_i \kappa(a_{i(-1)}b_{i(-1)})\ot \dif a_{i(0)}\ot
                 \dif b_{i(0)}}{s\ot t}=0 \quad
               \mbox{ for all $s\in T_\Omega$, $t\in T$}
     \end{align*}          
     implies
       $\ppair{\sum_i \omega_L(a_i^+b_{i(-1)})
              \ot \omega_L(b_{i(0)})}{s\ot t}=0$
        for all $s,t\in T_\Omega$,
     and hence
     \begin{align*}
       \sum_i \omega_L(a_i^+b_{i(-1)})\ot \omega_L(b_{i(0)})=0.
     \end{align*}
  \end{proof}
  
\begin{cor}\label{dimT0}
     Let $\B\subset \cA$ be as in Subsection \ref{qhs} and let
     $\Omega$, $\Gamma$, $T_\Omega$, $T$, and $T_0$ be as in
     Lemma \ref{ppair}. Assume that $\Gamma \B^+\subset \B^+\Gamma$.
     Then
     \begin{align*}
       \dim_\C T_0=\dim_\C T(\dim_\C T_\Omega-\dim_\C T).
     \end{align*}  
\end{cor}  
\begin{proof}
  By Lemma \ref{ppair} and the Remark \ref{catequivremarks}(i) one gets
  \begin{align*}
    \dim_\C(T_\Omega\ot T)/T_0
    =\dim_\C({}^{\co \cA}\cA\ot_\B\Gamma\ot_\B\Gamma)=
    \dim_\C(\Gamma\ot_\B\Gamma)/(\B^+\Gamma\ot_\B\Gamma).
  \end{align*}
  The inclusion $\Gamma\B^+\subset \B^+\Gamma$ implies that the canonical map
  \begin{align*}
    (\Gamma\ot_\B\Gamma)/(\B^+\Gamma\ot_\B\Gamma)\rightarrow
    \Gamma/\B^+\Gamma\ot_\C \Gamma/\B^+\Gamma  
  \end{align*}
  is an isomorphism. Therefore $\dim_\C T_0=(\dim_\C T_\Omega)(\dim_\C T)
                          -(\dim \Gamma)^2$.
\end{proof}

\begin{cor}\label{dimGuw2}
   Let $\B\subset \cA$ be as in Subsection \ref{qhs} and let
     $\Omega$, $\Gamma$, $T_\Omega$, $T$, and $T_0$ be as in
     Lemma \ref{ppair}. Then the pairing (\ref{nondegppairing}) induces a
     pairing
     \begin{align}\label{dimGuw2pairing}
       (\cA\ot_\B\Gamma_u^{\wedge 2})\times\Lin_\C\left\{\sum_i s_i\ot t_i\in
     T_\Omega\ot T\,\bigg|\,\sum_i s_it_i\in T\right\}\bigg/ T_0\rightarrow \C
     \end{align}
     which is nondegenerate when restricted to
     ${}^{\co \cA}(\cA\ot_\B\Gamma_u^{\wedge 2})$
     in the first component. 
\end{cor}
\begin{proof}
  Note first that $\Gamma^{\wedge 2}_{\mathrm{u}}=\Gamma^{\ot 2}/\Lambda$ where
  $\Lambda\subset \Gamma^{\ot 2}$ is the \textit{left}
  $\B$-submodule
  generated
  by $\{\sum_i \dif a_i\ot\dif b_i\in \Gamma^{\ot 2}\,|\,
  \sum_i  a_i\dif b_i=0\}$. It suffices to show that with respect to
  the pairing (\ref{ppairingdef}) one has
  \begin{align*}
    (\cA\ot_\B \Lambda)^\perp=
    \left\{\sum_i s_i\ot t_i\in T_\Omega\ot T\,\bigg|\,
   \sum_i s_it_i\in T\right\}.
  \end{align*}
  Assume that $\sum_i a_i\dif b_i=0$. Then
  \begin{align*}
     \sum_{i,j}\ppair{1\ot \dif a_i\ot \dif b_i}{s_j\ot t_j}&=
     \sum_{i,j} s_j(a_i^+b_{i(-1)})t_j(b_{i(0)})\\
     &\overset{\small (\ref{adb=0})}{=}-\sum_{i,j} \vep(a_i)
     s_j(b_{i(-1)})t_j(b_{i(0)})\\
     &=-\sum_{i,j}  s_jt_j(\vep(a_i)b_i)=-\sum_{i,j}  s_jt_j(\vep(a_i)b_i^+).
  \end{align*}
  Hence (\ref{rid}) implies that
  $\sum_j s_j\ot t_j\in (\cA\ot_\B \Lambda)^\perp$
  if and only if $\sum_j s_j t_j\in T$.
\end{proof}  

\section{Differential Calculus on Quantized Irreducible Flag Manifolds}
\label{DCQIFM}
In the previous Section we have recalled basic notions and developed the
general theory necessary for the investigation of covariant DC on quantum
homogeneous spaces.
Now we turn to the concrete example of quantized flag manifolds.
We first collect some facts about the corresponding algebras.
Then using the tools from the previous section the canonical covariant
DC over irreducible quantized flag manifolds is constructed and investigated
in detail. 
\subsection{Quantized Flag Manifolds}\label{QFM}
\subsubsection{Homogeneous Coordinate Rings}\label{HomCoordRing}
The quantized homogeneous coordinate ring $\sqgp$ of a generalized flag
manifold $G/P_S$ is defined to be the subalgebra of $\cqg$ generated by the
matrix coefficients $\{c_{f,v_\lambda}^\lambda\,|\, f\in V(\lambda)^\ast\}$,
\cite{b-CP94}, \cite{a-LaksResh92}, \cite{a-TaftTo91}, \cite{a-Soib92},
where $v_\lambda$ is a highest weight vector of $V(\lambda)$.
As a $\uqg$-module algebra $\sqgp$ is isomorphic to 
$\bigoplus_{n=0}^\infty V(n\lambda)^\ast$, where
$\lambda=\sum_{s\notin S}\omega_s$, endowed with the Cartan multiplication
\begin{align*}
  V(n_1\lambda)^\ast \ot V(n_2\lambda)^\ast
  \rightarrow V((n_1+n_2)\lambda)^\ast.
\end{align*}
Recall that the subspace $V(2\lambda)\subset V(\lambda)\ot V(\lambda)$ is
the eigenspace of $\rh_{\lambda,\lambda}$ with corresponding eigenvalue
$q^{(\lambda,\lambda)}$.
It is known (\cite{a-TaftTo91}, \cite{a-Brav94}) that $\sqgp$ is quadratic,
more explicitly
\begin{align*}
  \sqgp\cong \C\big\langle f_1,\dots,f_N\big\rangle \big/
         \big(\sum_{k,l\in I}\rh^{ij}_{kl}f_kf_l-q^{(\lambda,\lambda)}
         f_if_j\big)
\end{align*}
where $\rh$, $N$, $I$ are as in Section \ref{subsRMatrices}.
Similarly the dual quantized homogeneous coordinate ring $\sqgpop$ of
$G/P_S$ is defined to be the subalgebra of $\cqg$ generated by
$\{c_{v,f_{-\lambda}}^{-w_0\lambda}\,|\, v\in V(\lambda)\}$ where
$f_{-\lambda}\in V(-w_0\lambda)\cong V(\lambda)^\ast$ denotes the lowest
weight vector dual to $v_\lambda$.  
In terms of generators and relations one has
\begin{align*}
  \sqgpop\cong \C\big\langle v_1,\dots,v_N\big\rangle \big/
         \big(\sum_{k,l\in I}\rc^{ij}_{kl}v_kv_l
         -q^{(\lambda,\lambda)}v_iv_j\big).
\end{align*}
The $\uqg$-module structure of $\sqgp$ and $\sqgpop$ is given by identifying
the generators $\{f_i\,|\, i\in I\}$ and $\{v_i\,|\, i\in I\}$ with the
bases of $V(\lambda)^\ast$ and $V(\lambda)$ chosen in Section
\ref{subsRMatrices}. For notational reasons in what follows suppose that
$v_\lambda=v_N$ and $f_\lambda=f_N$.

\subsubsection{The Subalgebra $\sqgpcc\subset \cqg$}\label{subalgebra}
The tensor product $\sqgp_\C:=\sqgp\ot\sqgpop$ can be endowed with a
$\uqg$-module algebra structure by
\begin{align}\label{vf}
v_i f_j:=q^{(\lambda,\lambda)} \sum_{k,l\in I}\rgm ^{ij}_{kl} f_k\ot v_l.
\end{align}
To simplify notation the tensor product symbol will be omitted in the
following. The algebra $\sqgp_\C$ admits a character $\vep$ defined by
$\vep(v_i)=\vep(f_i)=\delta_{iN}$.
Note that
\begin{align*}
 c:=\sum_{i\in I} v_if_i=q^{(\lambda,\lambda)}\sum_{i,k,l\in I}
        \rgm^{ii}_{kl}f_k  v_l
\end{align*}
is a central invariant element of $\sqgp_\C$ with $\vep(c)=1$. The quotient
$\sqgpcc:=\sqgp_\C/(c-1)$ is $\Z$-graded by $\deg f_i=1$, $\deg v_i=-1$.
Let $S_q^n[G/P_S]_\C^{c=1}\subset \sqgpcc$ denote the homogeneous component
of degree $n$ with respect to this grading. 
\begin{lemma}\label{sqgpccLemma}
  The $\uqg$-module algebra $\sqgpcc$ is isomorphic to the $\uqg$-module
  subalgebra of $\cqg$ generated by the matrix coefficients
  $c^\lambda_{f,v_N}$, $c^{-w_0\lambda}_{v,f_N}$, $f\in
  V(\lambda)^\ast$, $v\in V(\lambda)$. The isomorphism is given by
  \begin{align*}
     f\mapsto c_{f,v_N}^\lambda \qquad v\mapsto
     c^{-w_0\lambda}_{v,f_N}.
  \end{align*}  
\end{lemma}  
\begin{proof}
  The torus $\C[K_i,K_i^{-1}\,|\,i=1,\dots,r]\subset \uqg$ acts on $\cqg$ by
  \begin{align*}
    K_i\lact c_{f,v}^\lambda=c_{f,K_iv}^\lambda.
  \end{align*}
  The eigenspace decomposition with respect to this action induces a
  $\Z$-grading on the subalgebra $A\subset \cqg$ generated by the matrix
  coefficients $c_{f,v_N}^\lambda$, $c_{v,f_N}^{-w_0\lambda}$,
  $f\in V(\lambda)^\ast$, $v\in V(\lambda)$. More precisely,
  $A=\bigoplus_{n\in\Z}A_n$ where
  \begin{align*}
    A_n=\Lin_\C\{c_{f,v_N^{\ot k}}^{k\lambda}
    c_{v,f_N^{\ot l}}^{-lw_0\lambda}\,|\, k,l\ge 0, k-l=n\}.
  \end{align*}
  Note that $\Lin_\C\{c_{f,v_{\mu}}^{\mu}
    c_{g,v_{w_0\nu}}^{\nu}\,|\, f\in V(\mu)^\ast, g\in V(\nu)^\ast\}\cong
    V(\mu)^\ast \ot V(\nu)^\ast$
  where $v_\mu$ and $v_{w_0\nu}$ denote a highest weight vector of
  $V(\mu)$ and a lowest weight vector of $V(\nu)^\ast$, respectively.    
  Therefore the relation
  $q^{(\lambda,\lambda)}\sum_{i,k,l\in I}(\rg^-)_{kl}^{ii}c_{f_k,v_N}^\lambda
  c_{v_l,f_N}^{-w_0\lambda}=1$ implies that $A_n$ can be written as
  a direct limit
  \begin{align*}
    A_n\cong\lim_{k\rightarrow \infty} V(k\lambda)^\ast\ot V((k-n)\lambda).
  \end{align*}
  of vector spaces.
  By construction the homogeneous components $S_q^n[G/P_S]^{c=1}_\C$ of
  \begin{align*}
    \sqgpcc=\bigoplus_{n=0}^\infty S_q^n[G/P_S]^{c=1}_\C
  \end{align*}
  allow the same presentation.
\end{proof}

\subsubsection{Quantized Flag Manifolds in Terms of Generators and Relations}
Lemma \ref{sqgpccLemma} implies that the subalgebra
$S_q^0[G/P_S]^{c=1}_\C$ is isomorphic
to the subalgebra $\A_\lambda^q\subset \cqg$ generated by the elements
$z_{ij}:=c_{f_i,v_N}^\lambda c_{v_j,f_N}^{-w_0\lambda}$.
It follows from
\begin{align}\label{c-lamclam}
  c_{v_i,f_N}^{-w_0\lambda}c_{f_j,v_N}^\lambda=
  q^{(\lambda,\lambda)}\sum_{kl\in I}(\rg^-)_{kl}^{ij}c_{f_k,v_N}^\lambda
  c_{v_l,f_N}^{-w_0\lambda}
\end{align}
that the following relations hold in $\A_\lambda^q$:
\begin{align}\label{zrelationen}
\ph_{12}\ra_{23}zz=0,\quad \pc_{34}\ra_{23}zz=0,\quad
q^{(\lambda,\lambda)}\sum_{i,j\in I}C_{ij}z_{ij}&=1,
\end{align}  
where $\ph:=(\rh-q^{(\lambda,\lambda)}\id)$,
$\pc:=(\rc-q^{(\lambda,\lambda)}\id)$,
$C_{kl}:=\sum_{i\in I}(\rg^-)_{kl}^{ii}$, and leg-notation is applied in the
first two formulae. Thus, explicitly written the first two equations of
(\ref{zrelationen}) take the form
\begin{align*}
  \sum_{m,n,p,t\in I}\ph^{ij}_{nm}\ra^{mk}_{pt}z_{np}z_{tl}=0,\qquad
  \sum_{m,n,p,t\in I}\pc^{kl}_{mt}\ra^{jm}_{np}z_{in}z_{pt}=0.
\end{align*}  
Let $\tilde{\A}_\lambda^q$ denote the free algebra $\C\langle z_{ij}\rangle$
divided by the ideal generated by the relations
(\ref{zrelationen}). It follows from the Yang-Baxter-Equation
that the left $\uqg$-module homomorphisms
\begin{align*}
  &V(n\lambda)^\ast\ot V(n\lambda)\rightarrow \tilde{\A}_\lambda^q,\quad
  f_{i_1}\dots f_{i_n}v_{j_1}\dots v_{j_n}\mapsto\\
  &\quad q^{-n(n-1)(\lambda,\lambda)/2}
  \sum_{\makebox[0cm]{$k_1,\dots,k_n\atop l_1,\dots,l_n$}}
    \left(\prod_{k=1}^{n-1}\prod_{l=1}^{n-k}
  \ra_{n+k-l,n+k-l+1} \right)^{i_1\dots i_n j_q\dots j_n}
  _{k_1 l_1\dots k_n l_n}z_{k_1l_1}\dots z_{k_nl_n}
\end{align*}  
\begin{align*}
  &V(n\lambda)^\ast\ot V(n\lambda)\hookrightarrow
  V((n+1)\lambda)^\ast\ot V((n+1)\lambda),\quad
  f_{i_1}\dots f_{i_n}v_{j_1}\dots v_{j_n}\mapsto\\
  &\qquad q^{(n+1)(\lambda,\lambda)}
  \sum_{\makebox[0cm]{$k_1,\dots,k_n\atop i,l,m$}}
    \left(\prod_{k=1}^{n+1}
  \rgm_{k,k+1} \right)^{iii_1\dots i_n}_{m k_1 \dots k_n l}f_mf_{k_1}\dots
  f_{k_n}v_lv_{j_1}\dots v_{j_n}
\end{align*}  
are well defined. Thus one obtains a surjection 
\begin{align*}
  A_0\cong\underrightarrow{\lim} V(n\lambda)^\ast\ot V(n\lambda)\rightarrow
  \tilde{\A}_\lambda^q.
\end{align*}
Note that since $A_0\subset \cqg$ the isotypical components of the
$\uqg$-module $\underrightarrow{\lim} V(n\lambda)^\ast\ot V(n\lambda)$
are finite dimensional. 
As the homomorphism $ \tilde{\A}_\lambda^q\rightarrow \A^q_\lambda$,
$z_{ij}\mapsto z_{ij}$ is surjective and
$\A^q_\lambda\cong\underrightarrow{\lim} V(n\lambda)^\ast\ot V(n\lambda)$
this yields  $\tilde{\A}_\lambda^q\cong\A^q_\lambda$.
Define
\begin{align}\label{cqgldef}
  \C_q[G/L_S]=\{a\in \C_q[G]\,|\, a_{(1)}\,a_{(2)}(k)=\vep(k)a\quad\forall
  k\in K\},
\end{align}  
where $K:=\U_q(\lfrak_S)$ is the Hopf subalgebra of $\uqg$ generated by
the elements $\{K_i,K_i^{-1},E_j,F_j\,|\, j\in S,\,i=1,\dots,r\}$.
By construction $\C_q[G/L_S]$ is a left $\C_q[G]$-comodule algebra containing
$\A_\lambda^q$. The following proposition was proved in \cite{a-stok02p},
\cite{a-heko03p}.
\begin{proposition}
  $\A_\lambda^q\cong \C_q[G/L_S]$ as left $\C_q[G]$-comodule algebras.
\end{proposition}

\subsection{First Order Differential Calculus over $\cqgl$}\label{QFFODC}

\subsubsection{Notations and Conventions}\label{subsNotations}
From now on we assume that $G/P_S$ is an irreducible flag manifold, in
particular $S=\{\alpha_s\}$, $\lambda=\omega_s$ for a fixed
$s\in\{1,\dots,r\}$.
To simplify notation, the isomorphic $\cqg$-comodule algebras $\cqgl$,
$\A_\lambda^q$, and $\tilde{\A}_\lambda^q$ will be denoted by $\B$.
By definition of $\cqgl$ the algebra $\B$ is a quantum homogeneous spaces
in the sense of \ref{qhs}.
Further, define
$I_{(1)}:=\{i\in I\,|\,(\omega_s,\omega_s-\alpha_s-\wght(v_i))=0\}$.
Note that the elements of $I_{(1)}$ label a basis of the $\uqls$-submodules
$V(\omega_s)_{(1)}:=\Lin_\C\{v\in V(\omega_s)\,|\,
       \wght(v)=\omega_s-\alpha_s-\beta,\,\beta\in
       Q,(\omega_s,\beta)=0\}\subset V(\omega_s)$
and $V(\omega_s)_{(1)}^\ast:=\Lin_\C\{f\in V(\omega_s)^\ast\,|\,
       \wght(f)=-\omega_s+\alpha_s+\beta,\,\beta\in Q,(\omega_s,\beta)=0\}
       \subset V(\omega_s)^\ast.$
As in \ref{HomCoordRing} assume that $v_N=v_{\omega_s}$ is the highest weight
vector of $V(\omega_s)$ in the basis $\{v_i\,|\,i\in I\}$ 
chosen in \ref{subsRMatrices}.

Fix a reduced decomposition of the longest element of the Weyl group.
Let $E_\beta, F_\beta$, $\beta\in R^+$, denote the corresponding root
vectors in $U$ \cite[8.1]{b-CP94}, \cite[6.2]{b-KS}.

Finally we introduce the abbreviation $M:=\dim_\C \gfrak/\pfrak_S=\# \gproots$.

\subsubsection{FODC over $\sqgp$}       
In analogy to the construction of $\B$ via $\sqgp$ one can obtain
covariant FODC over $\B$ by first constructing covariant FODC over $\sqgp$.
Consider the left $\sqgp$-module $\Gamma_+$ generated by elements
$\dif f_i$, $i\in I$, and relations
\begin{align*}
 \sum_{i,j\in I}a_{ij}f_i\dif f_j=0\quad\mbox{if }
 \sum_{i,j\in I}a_{ij}f_i{\ot}  f_j\in V(\mu)^\ast{\subset}
  V({\omega_s})^\ast {\ot}
  V({\omega_s})^\ast,\,\mu\neq 2{\omega_s}, 2{\omega_s}{-}\alpha_s.
\end{align*}  
As the subspaces $V(2{\omega_s})$ and $V(2{\omega_s}-\alpha_s)$
of $V({\omega_s})\ot V({\omega_s})$ are uniquely determined by the respective
eigenvalues  $q^{({\omega_s},{\omega_s})}$ and 
$-q^{({\omega_s},{\omega_s})-(\alpha_s,\alpha_s)}$ of
$\rh_{{\omega_s},{\omega_s}}$
these relations are equivalent to
\begin{align}\label{fdfrels}
 \sum_{i,j\in I}\left[\ph \qh \right]^{kl}_{ij}f_i\dif f_j=0\qquad \forall
 k,l\in I
\end{align}  
where as above $\ph=(\rh-q^{({\omega_s},{\omega_s})}\id)$ and
$\qh:=(\rh+q^{({\omega_s},{\omega_s})-(\alpha_s,\alpha_s)}\id)$.
The left module $\Gamma_+$ can be endowed with an $\sqgp$-bimodule structure
by
\begin{align}\label{dff}
  \dif f_i f_j
  =q^{(\alpha_s,\alpha_s)-({\omega_s},{\omega_s})}\sum_{k,l\in I}\rh^{ij}_{kl}f_k
  \dif f_l.
\end{align}  
Indeed, it follows from the Yang-Baxter-Equation for $\rh$ that this right
module structure is well defined on $\Gamma_+$.
As $\Lin_\C\{\dif f_i\}\cong V({\omega_s})^\ast$ the bimodule $\Gamma_+$
inherits a $\uqg$-module structure.

Define a linear map
\begin{align*}
  \dif:\sqgp\rightarrow \Gamma_+
\end{align*}  
by $\dif (f_i):=\dif f_i$ and $\dif (ab)=\dif a\,b+a\dif b$ for all $a,b\in
\sqgp$. To verify that $\dif$ is well defined note first that for
$\sum_{i,j\in I}a_{ij}f_i\ot  f_j\in V(\mu)^\ast\subset V({\omega_s})^\ast \ot
  V({\omega_s})^\ast$, $\mu\neq 2{\omega_s}, 2{\omega_s}-\alpha_s$, one has
$\sum_{i,j\in I}a_{ij}f_i \dif f_j=\sum_{i,j\in I}a_{ij}\dif f_i\, f_j= 0$.
Moreover, if $\sum_{i,j\in I}a_{ij}f_i\ot
f_j\in V(2{\omega_s}-\alpha_s)^\ast\subset
V({\omega_s})^\ast \ot V({\omega_s})^\ast$ then
\begin{align*}
\sum_{i,j,k,l\in I}a_{ij}\qh^{ij}_{kl}f_k\ot f_l=0
\end{align*}
and therefore
\begin{align*}
\sum_{i,j\in I}a_{ij}\dif f_i\, f_j+ a_{ij}f_i\dif f_j=
\sum_{i,j,k,l\in I}a_{ij}\left(q^{(\alpha_s,\alpha_s)-({\omega_s},{\omega_s})}\rh+\id
  \right)^{ij}_{kl}f_k\dif f_l=0.
\end{align*}

\subsubsection{FODC over $\sqgp_\C$}  
One can use $\Gamma_+$ to construct a covariant FODC
$(\Gamma_{+,\C},\del)$ over
$\sqgp_\C$ as follows. The left $\sqgp_\C$-module
\begin{align*}
  \Gamma_{+,\C}:=\sqgp_\C\ot_{\sqgp} \Gamma_+\cong \sqgpop\ot_\C\Gamma_+
\end{align*}
can be endowed with a right $\sqgp_\C$ module structure by
\begin{align}\label{dfv}
  \dif f_i v_j =q^{-({\omega_s},{\omega_s})} \sum_{k,l\in I}\ra ^{ij}_{kl} v_k
  \dif f_l.
\end{align}  
The differential $\del:\sqgp_\C\rightarrow \Gamma_{+,\C}$ defined by
\begin{align*}
  \del(v_i)=0,\quad \del (f_i)=\dif f_i
\end{align*}  
and Leibniz rule is well defined in view of (\ref{vf}), (\ref{dfv}).

There exists a pairing
\begin{align}\label{gam+cpairing}
  &\langle\cdot,\cdot\rangle:\Gamma_{+,\C}\ot V(\omega_s)_{(1)}\rightarrow
  \C,\\
  &\langle w\dif f,v\rangle :=\vep(w) f(v)\quad
           \mbox{ for $w\in \sqgp_\C$, $f\in \sqgp$}\nonumber
\end{align} 
where $\vep$ and $V(\omega_s)_{(1)}$ have been defined in \ref{subalgebra}
and \ref{subsNotations}, respectively.
To verify that $\langle\cdot,\cdot\rangle$ is well defined note that $\langle
       f_i\dif f_j,v\rangle\neq 0$ implies $\wght(f_i)=-{\omega_s}$,
$\wght(f_j)=-{\omega_s}+\alpha_s+\beta$, $(\beta,\omega_s)=0$, but then
$f_i\ot f_j \in V(2{\omega_s})^\ast\oplus V(2{\omega_s}-\alpha_s)^\ast\subset
V({\omega_s})^\ast\ot V({\omega_s})^\ast$.

Equations (\ref{vf}), (\ref{fdfrels}), (\ref{dff}), and (\ref{dfv}) imply
\begin{align}\label{delcrels}
\del c\,v_k=v_k\del c,\quad \del c\,f_k=q^{(\alpha_s,\alpha_s)}f_k\del c,
\quad \dif f_k \, c=c\dif f_k +(q^{(\alpha_s,\alpha_s)}-1)f_k\del c.
\end{align}

\subsubsection{The FODC $\Gamma_\del$}
Let $\Lambda\subset \Gamma_{+,\C}$ denote the subbimodule generated by
$\del c$,
$(c-1)\Gamma_{+,\C}$, and $\Gamma_{+,\C}(c-1)$. Then 
$\Gamma_{+,\C}/\Lambda$ is a covariant FODC over $\sqgp_\C^{c=1}$ which by
(\ref{delcrels}) as a
left module is generated by $\dif f_i$, $i\in I$, and relations
(\ref{fdfrels}) and $\del c=0$.
As $\vep(v_i)\neq 0$ if and only if $v_i=v_N$ and $f_N(v)=0$ for all $v\in
V(\omega_s)_{(1)}$ one obtains $\langle \del c,v\rangle=0$ for all
$v\in V(\omega_s)_{(1)}$. Therefore the pairing
\begin{align}\label{gam/Mpairing}
  \langle\,,\,\rangle:\Gamma_{+,\C}/\Lambda\ot V(\omega_s)_{(1)}\rightarrow \C
\end{align}  
induced by (\ref{gam+cpairing}) is well defined.
Let $\Gamma_\del\subset \Gamma_{+,\C}/\Lambda$ denote the FODC over
$\B\subset\sqgpcc$ induced by $\Gamma_{+,\C}/\Lambda$.
\begin{proposition}\label{delProp}
\begin{itemize}
\item[(i)]
  As a left $\B$-module $\Gamma_\del$ is generated by the differentials
  $\del z_{ij}$, $i,j\in I$, and relations
  \begin{align}
  \ph_{12}\qh_{12}\ra_{23}z\del z&=0, \label{pqrel}\\
  \pc_{34}\ra_{23}z\del z&=0, \label{prel}\\
  \sum_{i,j\in I}C_{ij}\del z_{ij}&=0.\label{Cdelz}
  \end{align}
\item[(ii)]
  The right $\B$-module structure of $\Gamma_\del$ is given by
  \begin{align}\label{rightmod}
  \del z z = q^{(\alpha_s,\alpha_s)}
  \rg^-_{23}\rh_{12}\rc^-_{34}\ra_{23}z\del z.
  \end{align}  
\item[(iii)] $\dim \Gamma_\del=M$.
\item[(iv)] $\Gamma_\del \B^+=\B^+\Gamma_\del.$
\item[(v)] The quantum tangent space of $\Gamma_\del$ is
  $\Lin_\C\{F_\beta\,|\, \beta\in \overline{R^+_S}\}$.
\end{itemize}
\end{proposition}  
\begin{proof}
Note first that the relations (\ref{pqrel})-(\ref{rightmod}) hold by
construction in  the $\B$-bimodule $\Gamma_\del$. Thus, by (\ref{rightmod})
$\Gamma_\del$ is generated by $\{\del z_{ij}\,|\, i,j\in I\}$ as a left
$\B$-module. 
Moreover, by (\ref{dfv}) the restriction of
the pairing (\ref{gam/Mpairing}) to
$\Lin_\C\{\del z_{iN}\,|\,i\in I_{(1)}\}\times V(\omega_s)_{(1)}$
is nondegenerate. Therefore one obtains
 $\dim \Gamma_\del\ge \dim_\C V(\omega_s)_{(1)}=M$.

To prove (i)-(iii) consider the left $\B$-module $\Gamma_\del'$
generated by elements
$\del z_{ij}$, $i,j\in I$, and relations (\ref{pqrel})-(\ref{Cdelz}).
By the categorial equivalence Theorem \ref{catequiv} it suffices to verify that
$\dim_\C\Gamma_\del'/\B^+\Gamma_\del'\le M$.
To this end note first that (\ref{prel}) multiplied by $\rc^-_{34}$, the
relation $\vep(z_{ij})=\delta_{iN}\delta_{jN}$, and (\ref{rtabelle}) imply
\begin{align*}
  \sum_{m,n\in I}\ra^{jk}_{mn}z_{im}\del z_{nl}\in \B^+\Gamma_\del'\quad
  \mbox{for all } i,j,k,l\in I \mbox{ such that } l\neq N.
\end{align*}
In particular one obtains
\begin{align}\label{delzinB+Gam1}
  \del z_{kl}\in \B^+\Gamma_\del'\quad
  \mbox{for all } k,l\in I \mbox{ such that } l\neq N.
\end{align}
Moreover, one calculates
\begin{align}
q^{({\omega_s},{\omega_s})}C_{34}&\rg^-_{23}\rh_{12}\rc^-_{34}\ra_{23}z\del z
\nonumber\\
&\overset{\small (\ref{prel})}{=}
C_{34}\rg^-_{23}\rh_{12}\ra_{23}z\del z\nonumber\\
&\overset{\small (\ref{pqrel})}{=}
q^{2({\omega_s},{\omega_s})-(\alpha_s,\alpha_s)}
C_{34}\rg^-_{23}\rh^-_{12}\ra_{23} z\del z+
q^{({\omega_s},{\omega_s})}(1-q^{-(\alpha_s,\alpha_s)})
\underbrace{C_{34}z\del z}_{=0 \mbox{\scriptsize \,\, by (\ref{Cdelz})}}
\nonumber\\
&=q^{2({\omega_s},{\omega_s})-(\alpha_s,\alpha_s)}
C_{12}\rg^-_{23}\rc^-_{34}\ra_{23} z\del z\nonumber\\
&\overset{\small (\ref{prel})}{=}
q^{({\omega_s},{\omega_s})-(\alpha_s,\alpha_s)}C_{12}z\del z\nonumber\\
&\overset{\small (\ref{zrelationen})}{=}q^{-(\alpha_s,\alpha_s)}\del z.
\label{zdelz=delz}
\end{align}  
Here, the third equation follows from the relations
\begin{align}\label{Crels}
  C_{23}\rg^-_{12}=C_{12}\rc^-_{23},\qquad
  C_{23}\rh^-_{12}=C_{12}\rg^-_{23}
\end{align}
which hold as the braiding induced by the action of the universal
$R$-matrix is a natural isomorphism. In view of (\ref{braiding}) and
(\ref{rtabelle}) Equation
(\ref{zdelz=delz}) implies
\begin{align}\label{delzinB+Gam2}
(q^{2(\wght(v_i)-\wght(v_j),\omega_s)}-q^{-(\alpha_s,\alpha_s)})\del z_{ij}
  \in \B^+\Gamma_\del'\qquad \mbox{for all $i,j\in I$}.
\end{align}
The relations (\ref{delzinB+Gam1}) and (\ref{delzinB+Gam2}) lead to
\begin{align}\label{delzinB+Gam3}
  \del z_{ij}\in \B^+\Gamma_\del' \mbox{ if $j\neq N$ or $i\notin I_{(1)}$.}
\end{align}
This proves
$\dim \Gamma_\del'=\dim_\C\Gamma_\del'/\B^+\Gamma_\del'\le\dim_\C
V(\omega_s)_{(1)}=M$.

We now prove (iv). By the third relation of (\ref{zrelationen})
the ideal $\B^+\subset \B$ is generated by
$\{z_{ij}\,|\, i\neq N\mbox{ or } j\neq N\}$.
Equation (\ref{rightmod}) and (\ref{rtabelle}) imply
that $\del z_{ij}z_{kl}$ can be written as a linear combination of elements
$z_{mn}\del z_{pt}$ where $\wght (v_k)\pord \wght (v_m)$ and
$\wght (v_l)\pord \wght (v_n)$. This proves $\Gamma_\del
\B^+\subset \B^+\Gamma_\del$. The converse inclusion follows similarly from
$z\del z=q^{-(\alpha_s,\alpha_s)}
\rg^-_{23}\rh^-_{12}\rc_{34}\ra_{23}\del z z$.

To prove (v) let $T$ denote the quantum tangent space of
$\Gamma_\del$. Recall from \cite[Lemma 6]{a-HK-QHS} and Remark
\ref{catequivremarks}(i)  that the pairing
\begin{align*}
  \Gamma_\del/\B^+\Gamma_\del\times T\rightarrow \C,\quad
  (\dif b,f)\mapsto f(b)
\end{align*}
is nondegenerate.
Moreover, by \cite[Theorem 7.2]{a-heko03p} there exist precisely two
non-isomorphic covariant FODC of dimension $M$ over
$\B$. The corresponding quantum tangent spaces are
$T_\del=\Lin_\C\{F_\beta\,|\, \beta\in \overline{R^+_S}\}$ and
$T_{\delb}=\Lin_\C\{E_\beta\,|\, \beta\in \overline{R^+_S}\}$.
As $T_{\delb}$ vanishes on all $z_{iN}$, $i\in I_{(1)}$, relation
(\ref{delzinB+Gam3}) implies $T\neq T_{\delb}$.
\end{proof}  

\subsubsection{The FODC $\Gamma_\delb$}

There exists a second covariant FODC $\Gamma_{\delb}$ over $\B$ of dimension
$\dim \gfrak/\pfrak_S$. This calculus can be obtained from a covariant FODC
over $\sqgpop$ in the same way as $\Gamma_\del$ has been obtained from
$\Gamma_+$. 
In analogy to $\Gamma_+$ a left $\sqgpop$-module $\Gamma_-$ can be defined by
generators $\dif v_i$, $i\in I$, and relations
\begin{align}\label{vdvrels}
 \sum_{i,j\in I}\left[\pc \qc \right]^{kl}_{ij}v_i\dif v_j=0\qquad \forall
 k,l\in I,
\end{align}  
where as before $\pc=(\rc-q^{({\omega_s},{\omega_s})}\id)$ and
$\qc:=(\rc+q^{({\omega_s},{\omega_s})-(\alpha_s,\alpha_s)}\id)$.
The left module $\Gamma_-$ can be endowed with a $\sqgpop$-bimodule structure
by
\begin{align*}
  \dif v_i v_j
  =q^{({\omega_s},{\omega_s})-(\alpha_s,\alpha_s)}\sum_{k,l\in I}
  \rcm^{ij}_{kl}v_k\dif v_l.
\end{align*}
Defining the differential $\dif :\sqgpop\rightarrow \Gamma_-$
by $\dif(v_i)=\dif v_i$ and the Leibniz rule one obtains the desired covariant
FODC over $\sqgpop$. To construct a covariant FODC $(\Gamma_{-,\C},\delb)$
over $\sqgp_\C$ consider the left $\sqgp_\C$-module
$\Gamma_{-,\C}:=\sqgp_\C\ot_{\sqgpop} \Gamma_-\cong \sqgp\ot_\C\Gamma_-$.
Then $\Gamma_{-,\C}$ can be endowed with a right $\sqgp_\C$ module structure
by 
\begin{align}\label{dvf}
  \dif v_i f_j =q^{({\omega_s},{\omega_s})}
  \sum_{k,l\in I}\rgm ^{ij}_{kl}f_k\dif v_l.
\end{align}  
The differential $\delb:\sqgp_\C\rightarrow \Gamma_{-,\C}$ is defined by
\begin{align*}
  \delb(f_i)=0,\quad \del (v_i)=\dif v_i.
\end{align*}  
and Leibniz rule.
There exists a pairing 
\begin{align}\label{gam-cpairing}
  &\langle\,,\,\rangle:\Gamma_{-,\C}\ot V(\omega_s)_{(1)}^\ast\rightarrow \C,\\
  & \langle w\dif v,f\rangle :=\vep(w) f(v), \quad\mbox{ for } w\in \sqgp_\C,
  \, v\in \sqgpop.\nonumber
\end{align}  
In analogy to (\ref{delcrels}) one has
\begin{align}\label{delbcrels}
\delb c\,f_k{=}f_k\delb c,\quad
\delb c\,v_k{=}q^{-(\alpha_s,\alpha_s)}v_k\delb c,
\quad \dif v_k \, c{=}c\dif v_k +(q^{-(\alpha_s,\alpha_s)}{-}1)v_k\delb c.
\end{align}
Let $\Lambda\subset \Gamma_{-,\C}$ denote the subbimodule generated by
$\delb c$,
$(c-1)\Gamma_{-,\C}$, and $\Gamma_{-,\C}(c-1)$.
Then $\Gamma_{-,\C}/\Lambda$ is a
covariant FODC over $\sqgp^{c=1}_\C$ which as a left module is generated by
$\dif v_i$, $i\in I$, and relations (\ref{vdvrels}) and $\delb c=0$.
Again the pairing
\begin{align}
  \langle\,,\,\rangle:\Gamma_{-,\C}/\Lambda\ot V(\omega_s)_{(1)}^\ast
  \rightarrow \C
\end{align} 
induced by (\ref{gam-cpairing}) is well defined. Let $\Gamma_{\delb}\subset
\Gamma_{-,\C}/\Lambda$ denote the FODC over $\B$ induced by
$\Gamma_{-,\C}/\Lambda$.
\begin{proposition}\label{delbProp}
\begin{itemize}
\item[(i)]
  As a left $\B$-module $\Gamma_{\delb}$ is generated by the differentials
  $\delb z_{ij}$, $i,j\in I$, and relations
  \begin{align}
  \pc_{34}\qc_{34}\ra_{23}z\delb z&=0, \label{delbpqrel}\\
  \ph_{12}\ra_{23}z\delb z&=0, \label{delbprel}\\
  \sum_{i,j\in I}C_{ij}\delb z_{ij}&=0.\label{Cdelbz}
  \end{align}
\item[(ii)]
  The right $\B$-module structure of $\Gamma_{\delb}$ is given by
  \begin{align}\label{delbrightmod}
  \delb z z = q^{-(\alpha_s,\alpha_s)}
  \rg^-_{23}\rh_{12}\rc^-_{34}\ra_{23}z\delb z.
  \end{align}  
\item[(iii)] $\dim \Gamma_{\delb}=M$.
\item[(iv)] $\Gamma_{\delb} \B^+=\B^+\Gamma_{\delb}.$
\item[(v)] The quantum tangent space of $\Gamma_{\delb}$ is
  $\Lin_\C\{E_\beta\,|\, \beta\in \overline{R^+_S}\}$.
\end{itemize}
\end{proposition}

\begin{proof}
  The proof is performed in analogy to the proof of Proposition
  \ref{delProp}. The following remarks may be helpful.
  Let $\Gamma_{\delb}'$ denote the left $\B$-module generated by elements
  $\delb z_{ij}$, $i,j\in I$, and relations
  (\ref{delbpqrel})-(\ref{Cdelbz}). Then relation (\ref{delbprel}) implies
  $\delb z_{ij}=\sum_{k\in I} z_{ik}\delb z_{kj}$ and therefore
  \begin{align}\label{delbzinB+Gam}
    \delb z_{kl}\in \B^+\Gamma_{\delb}'\quad \mbox{ for all $k,l\in I$ such
      that $k\neq N$.}
  \end{align}
  Similarly to (\ref{zdelz=delz}) one calculates
  \begin{align*}
    q^{({\omega_s},{\omega_s})}
      C_{34}&\rg^-_{23}\rh_{12}\rc^-_{34}\ra_{23}z\delb z
     =q^{(\alpha_s,\alpha_s)}\delb z
  \end{align*}   
  which in view of (\ref{rtabelle}) implies
  \begin{align*}
    (q^{2(\wght(v_i)-\wght(v_j),{\omega_s})}
          -q^{(\alpha_s,\alpha_s)})\delb z_{ij}
  \in \B^+\Gamma_{\delb}'.
  \end{align*}
\end{proof}

\subsubsection{The FODC $\Gamma_\dif$}
To obtain a $q$-deformed analogue of classical K\"ahler differentials over
$\C[G/L_S]$ we consider the sum
\begin{align}\label{Gammad}
  \Gamma_\dif:=\Gamma_\del+\Gamma_{\delb}.
\end{align}  
\begin{cor}\label{d=del+delb}
  \begin{itemize}
    \item[(i)]  $\Gamma_\dif=\Gamma_\del\oplus \Gamma_{\delb}$, in particular
      as a left $\B$-module $\Gamma_\dif$ is generated by the elements $\del
      z_{ij},\delb z_{ij}$, $i,j\in I$. 
    \item[(ii)] $\dim \Gamma_\dif=2M$.
    \item[(iii)] $\Gamma_\dif\B^+=\B^+\Gamma_\dif$. 
  \end{itemize}  
\end{cor}
\begin{proof}
  Recall that a sum of covariant FODC over $\B$ is direct if and only if the
  sum of their quantum tangent spaces is direct in $\B^\circ$. Now all
  statements of the corollary follow from Proposition \ref{delProp} and
  \ref{delbProp}.
\end{proof}  
\subsection{Higher Order Differential Calculus}\label{HODC}
The aim of this subsection is to determine the dimensions of the homogeneous
components of the universal DC $\Gamma^\wedge_{\del,\mathrm{u}}$,
$\Gamma^\wedge_{\delb,\mathrm{u}}$, and $\Gamma^\wedge_{\dif,\mathrm{u}}$ corresponding to
the covariant FODC $\Gamma_\del$, $\Gamma_{\delb}$, and $\Gamma_\dif$,
respectively. As in the previous subsection we first focus on $\Gamma_\del$.

\subsubsection{The Differential Calculus $\Gamma^\wedge_{\del,\mathrm{u}}$}
Recall from Proposition \ref{delProp}(iv) that
$\B^+\Gamma_\del=\Gamma_\del \B^+$ and hence
$\Gamma^\wedge_{\del,\mathrm{u}}/\B^+\Gamma^\wedge_{\del,\mathrm{u}}$ is an algebra generated
by $V_\del:=\Gamma_\del/\B^+\Gamma_\del$. For $i\in \Ind_{(1)}$
let $x_i\in V_\del$
denote the equivalence class of $\del z_{iN}\in \Gamma_\del$.
Note that $V_\del$ is an irreducible $K$-module isomorphic to
$V(\omega_s)_{(1)}^\ast$
with one-dimensional weight spaces. Therefore each irreducible $K$-module
in $V_\del\ot V_\del$ occurs with multiplicity $\le 1$. Hence the following
notion makes sense. An irreducible $K$-submodule of $V_\del \ot V_\del $
is called \textit{(anti)symmetric} if the weight vectors of the corresponding
classical $U(\mathfrak{l}_S)$-module are (anti)symmetric. Let $V_\del \ot
V_\del =S_\del \oplus A_\del $ denote the decomposition into the symmetric
and antisymmetric subspaces.

In order to describe the algebra
$\Gduw /\B^+\Gduw $ in terms of generators
and relations it is useful to consider the $-\N_0 $-filtration $\cH $ on
the vector space $V(\omega_s)_{(1)}^*\ot V(\omega_s)_{(1)}^*$ defined by
\begin{align*}
  \cH_n(V(\omega_s)_{(1)}^*\ot V(\omega_s)_{(1)}^*)=\Lin_\C\{E_\beta f_N\ot
  E_\gamma f_N\,|\,\max(\hght(\beta),\hght(\gamma))\ge -n\}
\end{align*}
where $\hght(\sum_{i=1}^r n_i\alpha_i)=\sum_{i=1}^rn_i$.
Moreover, we introduce the following notation:
for any $\beta \in \Rqps $ set $x_\beta :=x_i$ where $i\in \Ind _{(1)}$ and
$\wght (f_i)=\wght (E_\beta f_N)$.

Consider the totally ordered abelian
semigroup
\begin{align*}
  \cN=\{(k,n_1,\dots,n_k)\,|\, k\in \N_0, n_i\in -\N, n_i\le n_j\,
\forall\, i<j\}
\end{align*}
with the lexicographic order. The sum of two elements of $\cN$ is defined by
\begin{align*}
  (k,n_1,\dots,n_k)+(l,m_1,\dots,m_l)=(k+l,r_1,\dots,r_{k+l})
\end{align*}  
where $r_1,\dots,r_{k+l}$ are the numbers $n_1,\dots,n_k,m_1,\dots,m_l$ in
increasing order. 
The filtration $ \cH$ induces an
$\cN$-filtration on the algebra $\Gduw/\B^+\Gduw$  defined by
\begin{align}
  \deg(x_{\gamma})=
  (1,-\hght(\gamma))\label{defH}.
\end{align}  
This filtration will also be denoted by
$\cH$. 

\begin{proposition}\label{hodelProp}
\begin{itemize}
  \item[(i)]
    The algebra $\Gduw $ is generated by the elements $z_{ij},\del z_{ij}$,
  $i,j\in \Ind $, and relations (\ref{zrelationen}),
  (\ref{pqrel})--(\ref{rightmod}), and
  \begin{align}\label{relGduw}
   \qh_{12}\ra_{23}\del z\wedge \del z=0,\quad
   \pc_{34}\ra_{23}\del z\wedge \del z=0.
  \end{align}
\item[(ii)] The algebra $\Gduw /\B ^+\Gduw $ is isomorphic to
$V_\del ^\ot /(S_\del )$.\\
\item[(iii)] In the associated graded algebra
  $\gr _{\cH }\Gduw /\B ^+\Gduw $ the
following relations hold:
\begin{align*}
  x_\beta \wedge x_\gamma +q^{(\beta ,\gamma )}x_\gamma \wedge x_\beta =0
  \quad \text{for all $\beta ,\gamma \in \Rqps $ s.~t.
    $\hght(\gamma)\le\hght(\beta)$.}
\end{align*}
\item[(iv)] The set $\{x_{i_1}\wedge x_{i_2}\wedge \cdots \wedge x_{i_k}\,|\,
i_1<i_2<\cdots <i_k\}$ is a basis of $\Gduw[k]/\B^+\Gduw[k]$.
In particular $\dim \Gduw[k]={M\choose k}$.
\end{itemize}
\end{proposition}

\begin{proof}
(i) Recall that by construction $\Gduw $ is the quotient of the tensor
algebra $\Gamma _\del ^\ot $ by the ideal generated by
\begin{align*}
  \left\{\sum _i\del a_i\ot \del b_i\,\bigg|\,\sum _ia_i\del b_i=0\right\}.
\end{align*}
By Proposition \ref{delProp}(i),(ii) this ideal is generated by
\begin{align*}
  \{\ph _{12}\qh _{12}\ra _{23}\del z\ot \del z,
    \pc _{34}\ra _{23}\del z\ot \del z,
    \del z\ot \del z+q^{(\alpha _s,\alpha _s)}
    \rg^-_{23}\rh_{12}\rc^-_{34}\ra_{23}\del z\ot \del z\}
\end{align*}
and therefore coincides with the ideal generated by
\begin{align*}
  \{\qh _{12}\ra _{23}\del z\ot \del z,\pc _{34}\ra _{23}\del z\ot \del z\}.
\end{align*}
(ii) We first prove that $\dim \Gduw[2]=M(M-1)/2$.
Let $T_\Omega ^\vep \subset \uqg $ denote the
right coideal generated by $\{K_\beta F_\beta \,|\,\beta \in \Rqps \}$.
Let $\Omega $
denote the left covariant FODC over $\cA $ with quantum tangent space
$T_\Omega =(T_\Omega ^\vep )^+$. By Proposition \ref{delProp}(v) the space
$T_\Omega |_\B=T_\del $ is the quantum tangent space of $\Gamma _\del $.
By Proposition \ref{induced} one has $\Omega |_\B =\Gamma _\del $ and
therefore Corollary \ref{dimGuw2} can be applied.

By Corollary \ref{dimT0} and Proposition \ref{delProp}(iv) one obtains
$\dim _\C T_0=M(\dim _\C T_\Omega -M)$. On the other hand consider the linear
map
\begin{align*}
  \mathrm{m}: T_\Omega \ot T_\del \to \Ubar _-/T_\del,\qquad s\ot t\mapsto st
\end{align*}
where $\Ubar_-=\uqg/\uqg(K^++\uqb^+)$.
If $\beta_1,\dots,\beta_{M}$ denote the elements of $\Rqps$ then by
\cite[Prop.~5.2]{a-heko03p} the map $\mathrm{m}$ satisfies
$\im (\mathrm{m})=\Lin _\C \{F_{\beta _i}F_{\beta _j}\,|\,i\le j\}$
and hence $\dim \,\im (\mathrm{m})=M(M+1)/2$. By Corollary \ref{dimGuw2} this
implies
\begin{align}\label{dimGduw2}
  \dim \,\Gduw[2]=M(M-1)/2.
\end{align}
By the first equation of (\ref{relGduw}) and (\ref{rtabelle}) the generators
$x_i$ of $\Gduw /\B^+\Gduw $ satisfy the relation
\begin{align*}
  \sum _{k,l\in \Ind _{(1)}}\qh ^{ij}_{kl}x_k\wedge x_l=0\qquad
  \text{for all $i,j\in \Ind $.}
\end{align*}
For $\cQ :=\Lin _\C \{\sum _{k,l\in \Ind _{(1)}}\qh ^{ij}_{kl}x_k\ot x_l
\,|\,i,j\in \Ind \}\subset V_\del\ot V_\del$ by (\ref{dimGduw2}) one has
$\dim _\C \cQ \le M(M+1)/2
=\dim _\C S_\del $. Moreover, both $S_\del $ and $\cQ $ are $K$-submodules
of $V_\del \ot V_\del $. Therefore it suffices to show that the dimension
of any weight space of $S_\del $ does not exceed the dimension of the
corresponding weight space of $\cQ $.

For any element $\sum _{i,j\in \Ind }a_{ij}f_i\ot f_j\in V(2{\omega_s} )^*
\subset V({\omega_s} )^*\ot V({\omega_s} )^*$ where $a_{ij}\in \C $ one has
\begin{align*}
  \sum _{i,j,k,l\in \Ind }a_{ij}\qh ^{ij}_{kl}f_k\ot f_l
  =q^{({\omega_s} ,{\omega_s} )}(1+q^{-(\alpha _s,\alpha _s)})\sum _{i,j\in \Ind }
  a_{ij}f_i\ot f_j
\end{align*}
and therefore $\sum _{i,j\in \Ind _{(1)}}a_{ij}x_i\ot x_j \in \cQ $.
For $\beta ,\gamma \in \Rqps $ one calculates
\begin{align}
 E_\gamma E_\beta (f_N\ot f_N)=&E_\gamma (q^{-(\alpha _s,\alpha _s)/2}
  E_\beta f_N\ot f_N+f_N\ot E_\beta f_N)\notag \\
 =&q^{-(\alpha _s,\alpha _s)}E_\gamma E_\beta f_N\ot f_N
 +q^{-(\alpha _s,\alpha _s)/2}E_\beta f_N\ot E_\gamma f_N\nonumber\\
 &+q^{(\beta ,\gamma )-(\alpha _s,\alpha _s)/2}E_\gamma f_N\ot E_\beta f_N
 +f_N\ot E_\gamma E_\beta f_N\notag\\
 &+\sum _{i,j=1}^na_{ij}E_{\beta _i}f_N\ot E_{\beta _j}f_N
\label{EEff}
\end{align}
where in the last term $\beta _i,\beta _j\in \Rqps $ such that
$\max(\hght (\beta _i),\hght(\beta_j))>\hght (\beta )$ and the complex
numbers $a_{ij}$ depend on
$\beta $ and $\gamma $. Then (\ref{EEff}) implies that for every
$\beta ,\gamma \in \Rqps $ with $\hght (\gamma )\le \hght (\beta )$ there
exists $v_{\beta ,\gamma }\in
  \cH_n(V(\omega_s)_{(1)}^*\ot V(\omega_s)_{(1)}^*)$,
$n<-\hght(\beta)$, such that
\begin{align}\label{xbxg}
  x_\beta \ot x_\gamma +q^{(\beta ,\gamma )}x_\gamma \ot x_\beta +v_{\beta
  ,\gamma }\in \cQ.
\end{align}
This implies that the dimension of any weight space of $S_\del$ does not
exceed the dimension of the corresponding weight space of $\cQ $.

(iii) follows immediately from (\ref{xbxg}).

(iv) By (ii) the assertion holds for $k\le 2$. Moreover, (iii) implies that
the set $\{x_{i_1}\wedge x_{i_2}\wedge \cdots \wedge x_{i_k}\,|\,
i_1<i_2<\cdots <i_k\}$ generates the vector spaces
$\Gamma^{\wedge k}_{\del,\mathrm{u}} /
\B^+ \Gamma^{\wedge k}_{\del,\mathrm{u}}$.
By the diamond lemma it suffices to prove the claim for $k=3$.
To this end define
$V_{\delb}:= \Gamma_\delb/\B^+\Gamma_\delb\cong V(\omega_s)_{(1)}$
and let $V_\delb\ot V_\delb=S_\delb\oplus A_\delb$ denote the decomposition
into the symmetric and antisymmetric subspaces. By \cite[Cor.~6.7]{a-heko03p}
the graded vector spaces
$V_{\delb}^\ot/(A_{\delb})$ and $\C[x_1,\dots,x_M]$ are isomorphic. Thus
one has
\begin{align*}
  \dim(A_\delb\ot V_\delb + V_\delb\ot A_\delb)=M^3-{M+2 \choose 3}.
\end{align*}
The canonical pairing between $V(\omega_s)_{(1)}^\ast$ and
$V(\omega_s)_{(1)}$ induces a nondegenerate pairing of $K$-modules
\begin{align*}
  V_\del^{\ot 3}\ot V_\delb^{\ot 3}\rightarrow \C,\quad
  (x'{\ot} x''{\ot} x''')\ot(y'''{\ot} y''{\ot} y')\mapsto
  x'(y')x''(y'')x'''(y''').
\end{align*}
With respect to this pairing the equation
\begin{align*}
  S_\del\ot V_\del\cap V_\del\ot S_\del=(A_\delb\ot V_\delb + V_\delb\ot
  A_\delb )^\perp
\end{align*}
holds. Therefore
\begin{align*}
  \dim_\C( S_\del\ot V_\del+ V_\del\ot S_\del)&=
  2M\dim_\C S_\del-\dim(S_\del\ot V_\del\cap V_\del\ot S_\del)\\
  &=M^2(M+1)-{M+2 \choose 3}={M \choose 3}
\end{align*}
which by (ii) implies the claim for $k=3$.
\end{proof}

\subsubsection{The Differential Calculus $\Gdbuw$}

The situation for $\Gamma_\delb$ is completely analogous.
By Proposition \ref{delbProp}(v) one has $\B^+\Gamma_\delb=\Gamma_\delb\B^+$
and hence $\Gdbuw/\B^+\Gdbuw$ is an algebra generated by
$V_\delb:=\Gamma_\delb/\B^+\Gamma_\delb$. For $i\in I_{(1)}$ let $y_i\in
V_\delb$ denote the equivalence class of $\delb z_{Ni}\in \Gamma_\delb$.
Moreover, we use the following notation:
for any $\beta \in \Rqps $ set $y_\beta :=y_i$ where $i\in \Ind _{(1)}$ and
$\wght (f_i)=\wght (E_\beta f_N)$.

As in the proof of Proposition \ref{hodelProp} let $V_\delb\ot
V_\delb=S_\delb\oplus A_\delb$ denote the decomposition into the symmetric
and antisymmetric $K$-submodules. The algebra $\Gdbuw/\B^+\Gdbuw$ can be
endowed
with an $\cN$-filtration defined by $\deg (y_\gamma)=(1,-\hght(\gamma))$.
The proof of the following Proposition is a word by word translation of the
proof of Proposition \ref{hodelProp}.

\begin{proposition}\label{hodelbProp}
  (i) The algebra $\Gdbuw $ is generated by the elements $z_{ij},\delb z_{ij}$,
  $i,j\in \Ind $, and relations (\ref{zrelationen}),
  (\ref{delbpqrel})--(\ref{delbrightmod}) and
  \begin{align}\label{relGdbuw}
   \ph_{12}\ra_{23}\delb z\wedge \delb z=0,\quad
   \qc_{34}\ra_{23}\delb z\wedge \delb z=0.
  \end{align}
(ii) The algebra $\Gdbuw /\B ^+\Gdbuw $ is isomorphic to
$V_\delb ^\ot /(S_\delb )$.\\
(iii) In the associated graded algebra $\gr _{\cH }\Gdbuw /\B ^+\Gdbuw $ the
following relations hold:
\begin{align*}
  y_\beta \wedge y_\gamma +q^{-(\beta ,\gamma )}y_\gamma \wedge y_\beta =0
  \quad \text{for all $\beta ,\gamma \in \Rqps $ s.~t.
    $\hght(\gamma)\le\hght(\beta)$.}
\end{align*}
(iv) The set $\{y_{i_1}\wedge y_{i_2}\wedge \cdots \wedge y_{i_k}\,|\,
i_1<i_2<\cdots <i_k\}$ is a basis of $\Gdbuw[k]/\B^+\Gdbuw[k]$.
In particular $\dim \Gdbuw[k]={M\choose k}$.
\end{proposition}

\subsubsection{Extending $\del$ and $\delb$ to $\Gdifuw$}
Our next aim is to obtain results for $\Gdifuw$ analogous to Propositions
\ref{delProp} and \ref{delbProp}.  To this end we first show that the
decomposition $\Gamma_\dif=\Gamma_\del\oplus \Gamma_\delb$ induces
differentials $\del$ and $\delb$ on $\Gdifuw$ such that $\dif=\del+\delb$ also
holds in higher degrees. 

\begin{proposition}\label{delex}
  There exists a uniquely determined linear map
  $\del:\Gdifuw\rightarrow \Gdifuw$
  such that
  \begin{itemize}
    \item[(i)] $\del(\B)\subset \Gamma_\del\subset\Gamma_\dif$ and
      $\del|_\B$ coincides with the differential $\del$ considered in
      Section \ref{QFFODC}.
    \item[(ii)] $\del(\dif a)=-\dif(\del a)$ for all $a\in \B$. 
    \item[(iii)] $(\Gdifuw=\bigoplus_{i\in \N_0}\Gdifuw[i],\del)$ is a
      differential graded algebra. 
  \end{itemize}
  The map $\del$ satisfies $\del\dif=-\dif\del$. 
\end{proposition}  

\begin{proof}
  Uniqueness holds as $\B$ and $\dif \B$ generate the algebra $\Gdifuw$.
  To prove existence the following auxiliary lemma is needed. Let
  $\TomE$ and $\TomF$ denote the intersection of $\ker \vep$ with the right
  coideal of $\uqg$ generated by $\{E_\beta\,|\,\beta\in \Rqps\}$
  and $\{K_\beta F_\beta\,|\,\beta\in \Rqps\}$, respectively. By
  (\ref{hopfstruc}) the sum $T_\Omega=\TomE+\TomF\subset\uqg$ is direct. 
 Moreover, $T_\Omega$ is the quantum tangent space of a left covariant FODC
 $\Omega$ over $\cA$ such that $\Omega|_\B=\Gamma_\dif$. Let
  $\pomE,\pomF:T_\Omega=\TomE\oplus\TomF\rightarrow T_\Omega$ and
  $\pi_E,\pi_F:T =T_\del\oplus T_\delb\rightarrow T$ denote the
  canonical projections onto $\TomE$, $\TomF$, $T_\del$, and $T_\delb$,
  respectively. Recall the pairings (\ref{GamTpairing}) and
  (\ref{ppairingdef}).
  \begin{lemma}\label{adjoint1}
    The pairings
    \begin{align*}
      \pair{\cdot}{\cdot}:\Gamma_\dif\times T\rightarrow \C,\qquad
      \ppair{\cdot}{\cdot}:\cA\ot_\B\Gamma_\dif\ot_\B \Gamma_\dif\times
      T_\Omega\ot T\rightarrow \C
    \end{align*}
    satisfy the relations
    \begin{gather}
      \pair{\del a}{t}=\pair{\dif a}{\pi_F t},\qquad
         \pair{\delb a}{t}=\pair{\dif a}{\pi_E t},\label{pairdel(b)at}\\
      \ppair{\del a\ot \rho}{s\ot t}=
                 \ppair{\dif a\ot\rho}{\pomF s\ot t},\label{ppair1}\\
      \ppair{\delb a\ot \rho}{s\ot t}=
                 \ppair{\dif a\ot\rho}{\pomE s\ot t},\label{ppair2}\\
      \ppair{\rho \ot\del a}{s\ot t}=
                 \ppair{\rho\ot \dif a}{s\ot \pi_F t}, \label{ppair3}\\
      \ppair{\rho \ot\delb a}{s\ot t}=
                 \ppair{\rho\ot \dif a}{s\ot \pi_E t} \label{ppair4} 
    \end{gather}
    for all $a\in \B$, $\rho \in\Gamma_\dif$, $s\in T_\Omega$, $t\in T$.
  \end{lemma}  
\begin{proof}[Proof of Lemma \ref{adjoint1}]
  To prove (\ref{pairdel(b)at}) recall that by (\ref{delzinB+Gam3}) and
  Proposition \ref{delProp}(i)
 \begin{align*}
   \del a\in \Lin_\C \{\del z_{iN}\,|\,i\in I_{(1)}\}+\B^+\Gamma_\del
   \quad\mbox{ for all $a\in \B$.}
 \end{align*}
 Moreover, (\ref{delbzinB+Gam}) implies $\delb z_{iN}\in \B^+\Gamma_\dif$ if
 $i\in I_{(1)}$. Therefore using Corollary \ref{d=del+delb}(i) and
 $\del z_{iN}=\dif z_{iN} -\delb z_{iN}$ one obtains
 \begin{align*}
   \del a\in \Lin_\C \{\dif z_{iN}\,|\,i\in I_{(1)}\}+\B^+\Gamma_\dif
   \quad\mbox{ for all $a\in \B$.}
 \end{align*}
 Since $(\pi_E t)(z_{iN})=0$ for all $i\in I_{(1)}$ this implies
 $\pair{\del a}{\pi_E t}=0$ and hence
 \begin{align*}
   \pair{\del a}{t}=\pair{\del a}{\pi_F t} \quad\mbox{ for all $a\in \B$,
   $t\in T$.}
 \end{align*}
 Analogously one obtains
 \begin{align*}
   \pair{\delb a}{t}=\pair{\delb a}{\pi_E t} \quad\mbox{ for all $a\in \B$,
   $t\in T$}
 \end{align*}
 which yields (\ref{pairdel(b)at}).
 The remaining formulae follow from (\ref{pairdel(b)at}) and the definition
 (\ref{ppairingdef}) using
 \begin{align*}\left.
   \begin{array}{c}
     \pomF s_{(0)}^+\ot s_{(1)}=(\pomF s)_{(0)}^+\ot(\pomF s)_{(1)},\\
   \pomE s_{(0)}^+\ot s_{(1)}=(\pomE s)_{(0)}^+\ot(\pomE s)_{(1)}
   \end{array}\right\}
     \quad\mbox{ for all $s\in T_\Omega$.}
 \end{align*}  
\end{proof}
  We continue with the proof of Proposition \ref{delex}. The first step
  to prove existence of the map $\del:\Gdifuw\rightarrow \Gdifuw$ is to show
  that
  \begin{align}\label{del1def}
    \del:\Gdifuw[1]\rightarrow \Gdifuw[2],\quad \del(a\dif b):=\del
    a\wedge\dif b-a\,\dif\del b
  \end{align}
  is well defined. Assume that $\sum_i a_i\dif b_i=0$. Then
  $\sum_i a_i\del b_i=0$ and hence
  \begin{align*}
     \sum_i (\del a_i\wedge\dif b_i-a_i\dif\del b_i)&=
     \sum_i (\del a_i\wedge\dif b_i+\dif a_i\wedge\del b_i)\\
     &= \sum_i (2 \del a_i\wedge\del b_i+\del a_i\wedge\delb b_i+
     \delb a_i\wedge\del b_i).
  \end{align*}  
  Observe that $\{\sum_i (2 \del a_i\wedge\del b_i+\del a_i\wedge\delb b_i+
     \delb a_i\wedge\del b_i) \,|\,\sum_i a_i\dif b_i=0\}$ is a left
     $\cA$-subcomodule of $\Gdifuw[2]$. Thus by Corollary \ref{dimGuw2}
     it suffices to show that
  \begin{align}\label{suffcond}
    \sum_{i,j}\ppair{(2 \del a_i\wedge\del b_i+\del a_i\wedge\delb b_i+
     \delb a_i\wedge\del b_i)}{s_j\ot t_j}=0
  \end{align}
  whenever $\sum_i a_i\dif b_i=0$ and $\sum_j s_j t_j\in T$.
  By Lemma \ref{adjoint1} the left hand side of (\ref{suffcond}) is equal to
   \begin{align*}
    \sum_{i,j}\ppair{ \dif a_i\wedge\dif b_i}{s_j\ot t_j+\pomF s_j\ot\pi_F t_j
      -\pomE s_j\ot \pi_E t_j }.
  \end{align*}
  As $\sum_j\pomF s_j \pi_F t_j$, $\sum_j\pomE s_j \pi_E t_j \in T$
  whenever $\sum_j s_j t_j\in T$, the relation
  $\sum_i \dif a_i\wedge \dif b_i=0$ implies (\ref{suffcond}) and therefore
  (\ref{del1def}) is well defined.

  To verify that $\del$ is well defined on $\Gdifuw$ it suffices to check
  that
  \begin{align*}
    \sum_i \del\dif a_i\wedge \dif b_i-\dif a_i\wedge\del\dif b_i=0
  \end{align*}
  whenever $\sum_ia_i\dif b_i=0$, $a_i,b_i\in \B$. This follows from
  \begin{align*}
    \del\dif a\wedge\dif b-\dif a\wedge \del\dif b=-\dif(\del a\wedge \dif b-
    a\dif \del b)=-\dif(\del(a\dif b))
  \end{align*}  
  for all $a,b\in \B$.

  The property $\del \dif \rho=-\dif \del\rho $ for all $\rho\in \Gdifuw[k]$
  is proved by induction over $k$. 
  
  Next we show that $(\del\circ\del) |_\B=0$. To this end we calculate the
  adjoint operators of $\dif$ and $\del$ with respect to the pairing
  (\ref{dimGuw2pairing}). Assume that $\sum_js_j t_j\in T$, $s_j\in
  T_\Omega$, $ t_j\in T$, then
  \begin{align*}
     \bigppair{\dif(a\dif b)}{\sum_j s_j\ot t_j}&=
     \bigppair{\dif a\ot\dif b}{\sum_j s_j\ot t_j}\\
     &\overset{\small (\ref{adbdcppair})}{=}
          \sum_j  s_j(a^+b_{(-1)})t_j(b_{(0)})\\
     &=\sum_j s_j(ab_{(-1)})t_j(b_{(0)})-\vep(a) s_jt_j(b)     
  \end{align*}
  and hence
  \begin{align}\label{difadjoint}
    \bigppair{\dif \rho}{\sum_j s_j\ot t_j}=\sum_j s_j(\rho_{(-1)})
         \pair{\rho_{(0)}}{t_j}-\bigpair{\rho}{\sum_j s_j t_j}
  \end{align}
  for all $\rho\in \Gdifuw[1]$. Similarly, for all $a,b\in \B$,
  using Lemma \ref{adjoint1}, one obtains
  \begin{align*}
    \bigppair{\del(a\dif b)}{\sum_j s_j\ot t_j}=&
    \bigppair{\del a\wedge\dif b+\dif a\wedge\del b-\dif(a\del b)}
          {\sum_j s_j\ot t_j}\\
    \overset{\small (\ref{difadjoint})}{=}&
    \bigppair{\dif a\ot\dif b}{\sum_j (\pomF s_j\ot t_j+ s_j\ot\pi_F t_j)}\\
    &-\sum_j s_j(a_{(-1)}b_{(-1)})\pair{a_{(0)}\del b_{(0)}}{t_j}
                                        +\bigpair{a \del b}{\sum_j s_jt_j}\\
    \overset{\small (\ref{difadjoint})}{=}&
      \sum_j \pomF s_j(a_{(-1)}b_{(-1)})\pair{a_{(0)}\dif b_{(0)}}{t_j}\\
    & -\bigpair{a\dif b}{\sum_j[(\pomF s_j)t_j+ s_j(\pi_F t_j)-\pi_F(s_jt_j)]}.
  \end{align*}
This leads to
\begin{align*}
  \bigppair{\del^2 a}{&\sum_j s_j\ot t_j}=\sum_j\pomF s_j(a_{(-1)})
  \pair{\del a_{(0)}}{t_j}\\
       & -\bigpair{\del a}{\sum_j[s_jt_j-\pomE s_j\pi_E t_j
                      +\pomF s_J\pi_F t_j-\pi_F(s_jt_j)]}\\
  =& \sum_j\pomF s_j(a_{(-1)})  \pi_F t_j(a_{(0)})
     -\bigpair{\dif a}{\sum_j \pomF s_j\pi_F t_j}=0.
\end{align*}

It remains to prove $\del^2\rho =0$ for all $\rho\in \Gdifuw[k]$, $k\ge 1$,
which is obtained by induction over $k$. Assume $\del^2 \omega=0$ then using
$\del\dif =-\dif \del$ one gets
\begin{align*}
  \del^2(\dif a\wedge \omega)=\del(- \dif \del a \wedge \omega
    -\dif a\wedge \del \omega) =\dif(\del^2 a)\wedge \omega +\dif a\wedge
    \del^2 \omega=0.
\end{align*}  
\end{proof}  
\begin{remark}\upshape{
For the map $\delb:\Gdifuw\rightarrow \Gdifuw$ defined by $\delb:=\dif-\del$
one immediately obtains properties analogous to Proposition \ref{delex}.
Moreover, one verifies that $\del\delb+\delb\del=0$.
}
\end{remark}  

\subsubsection{The Differential Calculus $\Gdifuw$}
Now we are prepared to write the algebra $\Gdifuw$ in terms of generators
and relations and to calculate $\dim \Gdifuw[k]$ for all $k$.
By Propositions \ref{delProp}(v) and \ref{delbProp}(v) one has
$\B^+\Gamma_\dif=\Gamma_\dif\B^+$
and hence $\Gdifuw/\B^+\Gdifuw$ is an algebra generated by
$\Gamma_\dif/\B^+\Gamma_\dif=\Gamma_\del/\B^+\Gamma_\del\oplus
\Gamma_\delb/\B^+\Gamma_\delb$. Recall that for $\beta\in \Rqps$ we write
$x_\beta$ and $y_\beta$ to denote the equivalence class of $\del z_{iN}\in
\Gamma_\del$ and $\delb z_{Ni}\in \Gamma_\delb$ for suitable $i\in I_{(1)}$,
respectively.
The algebra $\Gdifuw/\B^+\Gdifuw$ can be endowed with an $\cN$-filtration
defined by $\deg (x_\gamma)=(1,-\hght(\gamma))=\deg(y_\gamma)$.
As before this $\cN$-filtration
will be denoted by $\cH$.
\begin{proposition}\label{hodProp}
  \begin{itemize}
    \item[(i)] The algebra $\Gdifuw$ is generated by the elements
    $z_{ij}$, $\del z_{ij}$, $\delb z_{ij}$, $i,j\in I$, and relations
    (\ref{zrelationen}), (\ref{pqrel}) -- (\ref{rightmod}),
    (\ref{delbpqrel}) -- (\ref{delbrightmod}), (\ref{relGduw}),
    (\ref{relGdbuw}), and
    \begin{align}\label{delbdelrel}
       \delb z {\wedge} \del z =
          -q^{-(\alpha_s,\alpha_s)}T^-_{1234}\del z {\wedge} \delb z +
       q^{(\omega_s,\omega_s)-(\alpha_s,\alpha_s)}z C_{12}T^-_{1234}\del z
          {\wedge}\delb z
    \end{align}
    where $T^-_{1234}=\rg^-_{23}\rh^-_{12}\rc_{34}\ra_{23}$.
    \item[(ii)] The algebra $\Gdifuw /\B ^+\Gdifuw $ is isomorphic to
    $(V_\del\oplus V_\delb) ^\ot /(S_\del+S_\delb+J )$ where
    $J\subset (V_\del\ot V_\delb)\oplus (V_\delb\ot V_\del)$ is the subspace
    spanned by all
    expressions of the form
    \begin{align}\label{yx=Rxy}
      y_i\ot x_j+q^{(\omega_s,\omega_s)-(\alpha_s,\alpha_s)}
      \sum_{k,l\in I_{(1)}} \rgm^{ij}_{kl}x_k\ot y_l, \quad i,j \in I_{(1)}.
    \end{align}  
    \item[(iii)] In the associated graded algebra
      $\gr _{\cH }\Gdifuw /\B ^+\Gdifuw $ the following relations hold:
    \begin{align*}
        y_\beta \wedge y_\gamma +q^{-(\beta ,\gamma )}
             y_\gamma \wedge y_\beta =0,\\
        x_\beta \wedge x_\gamma +q^{(\beta ,\gamma )}
             x_\gamma \wedge x_\beta =0
    \end{align*}         
    for all $\beta ,\gamma \in \Rqps $ such that
    $\hght(\gamma)\le\hght(\beta)$, and
     \begin{align*}
        y_\beta \wedge x_\gamma +q^{-(\beta ,\gamma )}
             x_\gamma \wedge y_\beta =0
      \end{align*}         
    for all $\beta ,\gamma \in \Rqps $.        
    \item[(iv)] For all $k\in \N_0$ the canonical map
      \begin{align}\label{isom}
         \bigoplus_{i+j=k}\Gduw[i]/\B^+\Gduw[i]\ot\Gdbuw[j]/\B^+\Gdbuw[j]
                       \rightarrow \Gdifuw[k]/\B^+\Gdifuw[k]
      \end{align}
      is an isomorphism. In particular $\dim \Gdifuw[k]={2M\choose k}$.
  \end{itemize}    
\end{proposition}  
\begin{proof}
 (i) By Corollary \ref{d=del+delb} the algebra $\Gdifuw$ is generated by the
  elements $z_{ij}$, $\del z_{ij}$, $\delb z_{ij}$, $i,j\in I$.
  Moreover, Propositions
  \ref{delProp} and \ref{delbProp} imply that the  relations
   (\ref{pqrel}) -- (\ref{rightmod}) and
   (\ref{delbpqrel}) -- (\ref{delbrightmod}) hold. Applying $\del$ and $\delb$
   one obtains (\ref{relGduw}) and (\ref{relGdbuw}). In the following we
   verify (\ref{delbdelrel}).

   Using (\ref{Crels}) and applying $C_{23}\rc^-_{34}$ to (\ref{prel}) and
   $C_{23}\rh^-_{12}$ to
   (\ref{delbprel}) one obtains for $D:=C\ra$ the relations
   \begin{align*}
     D_{23}z\,\del z&=0,&D_{23}z\,\delb z&=\delb z.
   \end{align*}  
   Leibniz rule for $\delb$ yields $D_{23}\delb z\,z=0$ and
   $D_{23}\delb z\wedge\del z=D_{23}z\,\del\delb z$.
   Thus one gets
   \begin{align*}
     \del\delb z=D_{23}\del(z\delb z)=D_{23}(\del
     z\wedge\delb z+\delb z\wedge \del z).
   \end{align*}
   With the abbreviation $T_{1234}=\rg^-_{23}\rh_{12}\rc^-_{34}\ra_{23}$
   using (\ref{rightmod}) and (\ref{delbrightmod}) this leads to
   \begin{align}
     \del\delb z\,z &=D_{23}(\del z\wedge\delb z+\delb z\wedge \del z)z
     =D_{23}T_{3456}T_{1234}(z\del z\wedge\delb z+z\delb z\wedge \del z)
       \nonumber\\
     &=D_{23}\rg^-_{45}\rh_{34}\rg^-_{23}\rh_{12}
             \rc^-_{56}\ra_{45}\rc^-_{34}\ra_{23}
     (z\del z\wedge\delb z+z\delb z\wedge \del z)\nonumber\\
     &=\rg^-_{23}D_{34}\rh_{12}\rc^-_{56}\ra_{45}\rc^-_{34}\ra_{23}
     (z\del z\wedge\delb z+z\delb z\wedge \del z)\nonumber\\
     &=\rg^-_{23}\rh_{12}\rc^-_{34}D_{45}\ra_{23}
     (z\del z\wedge\delb z+z\delb z\wedge \del z)\nonumber\\
     &=T_{1234}z\del\delb z \label{deldelbzz}
   \end{align}
   where the relations
   \begin{align*}
     D_{12}\rh_{23}\rg^-_{12}=D_{23},\qquad D_{12}\ra_{23}\rc^-_{12}=D_{23}
   \end{align*}
   have been used. Now $\delb$ is applied to (\ref{rightmod}) which leads to
   \begin{align*}
      \delb\del z\,z-\del z\wedge  \delb z=
             q^{(\alpha_s,\alpha_s)}T_{1234} \delb z\wedge \del z
            +q^{(\alpha_s,\alpha_s)}T_{1234}  z \delb\del z.
   \end{align*}
   In view of (\ref{deldelbzz}) multiplication by $T_{1234}^-$ yields
   \begin{align}\label{zdelbdelz}
     (1-q^{(\alpha_s,\alpha_s)})z\delb \del z= T_{1234}^- \del z\wedge\delb z+
     q^{(\alpha_s,\alpha_s)}\delb z \wedge\del z.
   \end{align}  
   Application of $C_{12}$ leads to
   \begin{align*}
     q^{-(\omega_s,\omega_s)}(1-q^{(\alpha_s,\alpha_s)})\delb \del z=
          C_{12}T_{1234}^-\del z\wedge\delb z.
   \end{align*}
   Inserting this formula in (\ref{zdelbdelz}) one finally gets the desired
   Equation (\ref{delbdelrel}).

   Now Propositions \ref{hodelProp}(iv), \ref{hodelbProp}(iv) and
   (\ref{delbdelrel}) imply
   \begin{align*}
     \dim \Gamma_\dif^{\ot 2}/\Lambda\le {2M \choose 2}
   \end{align*}
   where $\Lambda\subset\Gamma_\dif^{\ot 2}$ denotes the $\B$-bimodule
   corresponding
   to the relations (\ref{relGduw}), (\ref{relGdbuw}), and (\ref{delbdelrel}).
   Let $\Omega$ denote the left covariant FODC over $\cA$ defined in the
   proof of Proposition \ref{delex}. Recall that $\Omega |_\B=\Gamma_\dif$
   and $T_\Omega|_\B=T$ and therefore Corollary \ref{dimGuw2} can be
   applied.
   By Corollaries \ref{dimT0} and \ref{d=del+delb}(iii) one obtains
   $\dim _\C T_0=2M(\dim _\C T_\Omega -2M)$.
   On the other hand let $\beta_1,\dots,\beta_M$ denote the elements of
   $\Rqps$. By \cite[Prop.~5.2]{a-heko03p} the linear map
\begin{align*}
  \mathrm{m}: T_\Omega \ot T \to \Ubar/T,\qquad s\ot t\mapsto st.
\end{align*}
satisfies
$\im (\mathrm{m})=\Lin _\C \{F_{\beta _i}F_{\beta _j},
E_\beta F_\gamma,E_{\beta _i}E_{\beta _j}\,|\,i\le j,\,\beta,\gamma\in \Rqps\}$
and therefore $\dim \,\im (\mathrm{m})=2M(2M+1)/2$. By Corollary \ref{dimGuw2} 
one obtains
\begin{align}\label{dimGam2}
  \dim \,\Gdifuw[2]={2M \choose 2}.
\end{align}
This implies $\Gamma^{\ot 2}_\dif/\Lambda=\Gdifuw[2]$ and completes the proof
of (i) as $\Gdifuw$ is a quadratic algebra. 

(ii)  The algebra $\Gdifuw /\B ^+\Gdifuw $ is generated by the elements
$x_i,y_i$, $i=1,\dots,M$, and the relations induced by (\ref{relGduw}),
(\ref{relGdbuw}), and (\ref{delbdelrel}). It has already been stated in
Proposition \ref{hodelProp} and Proposition \ref{hodelbProp} that the
relations induced by (\ref{relGduw}) and (\ref{relGdbuw}) are obtained by
setting $S_\del\subset V_\del\ot V_\del$ and
$S_\delb\subset V_\delb\ot V_\delb$ equal to zero. On the other hand 
(\ref{rtabelle}) and (\ref{delbdelrel}) imply that (\ref{yx=Rxy})
vanishes in $\Gdifuw /\B ^+\Gdifuw $. By (\ref{dimGam2})
there can be no more quadratic relations.

(iii) The first two relations of (iii) have already been stated in Propositions
\ref{hodelProp} and \ref{hodelbProp}.
In view of (\ref{rtabelle}) and the definition of the filtration $\cH$ the
last relation follows from (\ref{yx=Rxy}).

(iv) By Equation (\ref{relGduw}), (\ref{relGdbuw}), and (\ref{delbdelrel})
the vector space $\Gdifuw/\B^+\Gdifuw$ is a quotient of the tensor product
$(\Gduw/\B^+\Gduw)\ot(\Gdbuw/\B^+\Gdbuw)$. Moreover, the map (\ref{isom})
is an isomorphism if and only if there exists a product $\wedge$ on 
$(\Gduw/\B^+\Gduw)\ot(\Gdbuw/\B^+\Gdbuw)$ which extends the algebra
structures of $\Gduw/\B^+\Gduw$ and $\Gdbuw/\B^+\Gdbuw$ and satisfies
(\ref{yx=Rxy}). Existence of the product $\wedge$ follows from Lemma
\ref{restrictedR} and the naturality of the braiding (\ref{braiding})
of $U_q([\lfrak_S,\lfrak_S])$.
\end{proof}
\begin{remark}\upshape{
(i) For $\gfrak= \mathfrak{o}_n$, $s=1$, and $\gfrak=\mathfrak{sl}_n$ the
   algebras $\Gduw/\B^+\Gduw$ are well known examples of quantized
   exterior algebras \cite[Def.~7.4.4]{b-CP94}, \cite{a-FadResTak1}.

(ii) Propositions \ref{hodProp}(iv), \ref{hodelProp}(iv), and
\ref{hodelbProp}(iv) imply that $\Gdifuw[2M]/\B^+\Gdifuw[2M]$ is
a one dimensional trivial $K$-module. Thus $\Gdifuw[2M]$ is a free left
$\B$-module generated by one left coinvariant element.
In contrast the covariant differential calculi $\Gduw$ and $\Gdbuw$ do not
admit a volume form.
}
\end{remark}

\bibliographystyle{amsalpha}
\bibliography{litbank2}

\end{document}